\documentclass[11pt,leqno,oneside,letterpaper]{amsart}

\usepackage[margin=1.15in]{geometry}
\usepackage{amsmath,amsfonts,amssymb,amscd,amsthm,amsbsy,epsf,graphics} 
\usepackage{stmaryrd}
\usepackage{setspace}
\usepackage{comment}
\usepackage{color}
\usepackage[colorlinks,linkcolor=blue,citecolor=blue]{hyperref} 
\usepackage{epsfig} 
\usepackage{tikz-cd}
\usepackage{bm}
\usepackage{graphicx}  
\usepackage[caption=false]{subfig}
\usepackage{lipsum}
\usepackage{enumitem} 

\usepackage[square,numbers]{natbib}

\makeatletter
\def\l@subsection{\@tocline{2}{0pt}{4pc}{5pc}{}}
\let\oldtocsection=\tocsection
\let\oldtocsubsection=\tocsubsection
\let\oldtocsubsubsection=\tocsubsubsection
 
\renewcommand{\tocsection}[2]{\hspace{0em}\oldtocsection{#1}{#2}}
\renewcommand{\tocsubsection}[2]{\hspace{0em}\oldtocsubsection{#1}{#2}}
\renewcommand{\tocsubsubsection}[2]{\hspace{2em}\oldtocsubsubsection{#1}{#2}}

\newtheorem{thm}{Theorem}[section]
\newtheorem*{thm*}{Theorem}
\newtheorem{cor}[thm]{Corollary}
\newtheorem{lemma}[thm]{Lemma}
\newtheorem{prop}[thm]{Proposition}
\newtheorem{conj}[thm]{Conjecture}

\newtheorem{innercustomthm}{Theorem}
\newenvironment{customthm}[1]
  {\renewcommand\theinnercustomthm{#1}\innercustomthm}
  {\endinnercustomthm}

\newtheorem{innercustomcor}{Corollary}
\newenvironment{customcor}[1]
{\renewcommand\theinnercustomcor{#1}\innercustomcor}
{\endinnercustomcor}

\theoremstyle{definition}
\newtheorem{question}[thm]{Question}

\newtheorem{remark}[thm]{Remark} 
\newtheorem{remarks}[thm]{Remarks}

\newtheorem*{acknowledgement*}{Acknowledgements}

\newtheoremstyle{cases}
  {12pt plus 6 pt}
  {2pt}
  {\bfseries}   
  {}
  {\bfseries}
  {.}
  {.5em}
  {}
\theoremstyle{cases}

\newtheorem{case}{Case}

\numberwithin{subcase}{case} 
\numberwithin{subsubcase}{subcase}
\numberwithin{equation}{subsection} 

\def\bQ{{\mathbb Q}}
\def\bZ{{\mathbb Z}}

\graphicspath{{figures/}} 

\begin{document}

\title[JSJ decompositions, Dehn surgery and the $L$-space conjecture]{JSJ decompositions of knot exteriors, Dehn surgery and the $L$-space conjecture}

\author[Steven Boyer]{Steven Boyer} 
\thanks{Steven Boyer was partially supported by NSERC grant RGPIN 9446-2008}
\address{D\'epartement de Math\'ematiques, Universit\'e du Qu\'ebec \`a Montr\'eal, 201 President Kennedy Avenue, Montr\'eal, Qc., Canada H2X 3Y7.}
\email{boyer.steven@uqam.ca}
\urladdr{http://www.cirget.uqam.ca/boyer/boyer.html}

\author[Cameron McA. Gordon]{Cameron McA. Gordon}
\address{Department of Mathematics, University of Texas at Austin, 1 University Station, Austin, TX 78712, USA.}
\email{gordon@math.utexas.edu}
\urladdr{}

\author[Ying Hu]{Ying Hu}
\address{Department of Mathematical and Statistical Sciences, University of Nebraska Omaha, 6001 Dodge Street, Omaha, NE 68182-0243, USA.}
\email{yinghu@unomaha.edu}
\urladdr{https://yinghu-math.github.io}

\thanks{2020 Mathematics Subject Classification.  Primary 57M05, 57M50, 57M99}

\thanks{Key words: JSJ-decomposition, Dehn surgery, satellite knots, $L$-space knots, left-orderable, taut foliations, $L$-space conjecture}

\begin{abstract}
In this article, we apply slope detection techniques to study properties of toroidal $3$-manifolds obtained by performing Dehn surgeries on satellite knots in the context of the $L$-space conjecture. We show that if $K$ is an $L$-space knot or admits an irreducible rational surgery with non-left-orderable fundamental group, then the JSJ graph of its exterior is a rooted interval. Consequently, any rational surgery on a composite knot has a left-orderable fundamental group. This is the left-orderable counterpart of Krcatovich's result on the primeness of $L$-space knots, which we reprove using our methods. Analogous results on the existence of co-orientable taut foliations are proved when the knot has a fibred companion. Our results suggest a new approach to establishing the counterpart of Krcatovich's result for surgeries with co-orientable taut foliations, on which partial results have been achieved by Delman and Roberts. Finally, we prove results on left-orderable $p/q$-surgeries on knots with $p$ small.  
\end{abstract}

\maketitle

\section{Introduction} 

We call a closed, connected, orientable, irreducible $3$-manifold $LO$ if its fundamental group is left-orderable \footnote{We adopt the usual convention in the context of the $L$-space conjecture that the trivial group is not left-orderable. See for instance \cite[\S 6.1]{BGH21}. Thus $S^3$ is not $LO$.}, $CTF$ if it admits a co-oriented taut foliation, and $NLS$ if it is not a Heegaard Floer $L$-space. The $L$-space conjecture contends that these three conditions are equivalent (\cite{BGW13}, \cite{Juhasz2015}). 

The notion of slope detection (\S \ref{subsec: slope detection}) was first introduced by Boyer and Clay in \cite{BC17} to study the $L$-space conjecture for graph manifolds (also see \cite{HRRW15}), and was further exploited by this paper's authors in \cite{BGH21} to study it in the context of toroidal $3$-manifolds. In particular the techniques were applied to toroidal rational homology spheres which have small first homology or which arise as cyclic branched covers of satellite links. A key aspect of this method is that it provides a unified approach to the verification of each of the three properties which arise in the conjecture.  

This paper illustrates the efficacy of the method through the study of Dehn surgery on satellite knots. Though the initial focus of this article  was to prove results on surgeries which yield $LO$-manifolds, the analogous behaviour of detection and gluing in the $LO, NLS$, and $CTF$ contexts (cf. \cite[\S 1]{BGH21}) means that the arguments used in the $LO$ case could be applied mutatis mutandis to prove results on $NLS$-surgery and $CTF$-surgery. This approach leads to new theorems on $LO$-, $NLS$-, and $CTF$-surgery, as well as new proofs of some known results on $NLS$-surgery. 

In stating our results we use the following conventions. A {\em knot} will mean a knot in the $3$-sphere. Given a knot $K$, $X(K)$ will denote its exterior and for $r \in \mathbb Q_\infty = \mathbb Q \cup \{\frac10\}$, we use $K(r)$ to denote the manifold obtained by $r$-surgery on $K$. When $r \in \mathbb Q$ is expressed as a reduced fraction $p/q$, we assume that $p \geq 0$ unless otherwise indicated. 

We first show that having an irreducible rational Dehn surgery that is either not $LO$ or an $L$-space puts strong restrictions on the essential tori in the exterior of a knot.

\begin{customthm}{\ref{thm: split}}
Suppose that the exterior $X(K)$ of a knot $K$ contains disjoint essential tori $T_1, T_2$ which together with $\partial X(K)$ cobound a connected submanifold of $X(K)$. Let $r \in \mathbb Q$ and assume that $r$ is not the cabling slope if $K$ is a cable knot. Then $K(r)$ is $LO$ and $NLS$.
\end{customthm}

Since cable knots are prime, the following corollary is immediate. See Figure \ref{fig: composite knots}.

\begin{customcor}{\ref{cor: nlo irred prime}}
All rational surgeries on a composite knot are $LO$ and $NLS$. 
\end{customcor}

To put this corollary in context, recall that an $L$-space knot is a knot which admits a rational $L$-space Dehn surgery. Corollary \ref{cor: nlo irred prime} then reproves a result of Krcatovich, which states  that $L$-space knots are prime \cite{Krc15}, and establishes the analogue for $LO$-surgery. 

The {\it JSJ graph} of (the exterior of) a knot $K$ is the rooted tree dual to the JSJ tori of $X(K)$, where the root vertex corresponds to the piece of the JSJ decomposition containing $\partial X(K)$. Any finite rooted tree can arise as the JSJ graph of a knot. We say that a JSJ graph is a {\it rooted interval} when it is homeomorphic to an interval with root at an endpoint. It is easy to see that the JSJ graph of $K$ is a rooted interval if and only if the exterior of $K$ does not contain essential tori $T_1$ and $T_2$ as in Theorem \ref{thm: split}. See Figure \ref{fig: jsj graph}.

Since an $L$-space knot has infinitely many $L$-space surgeries, we have the following immediate consequences of Theorem \ref{thm: split}.

\begin{customthm}{\ref{thm: jsj graph interval $L$-space knot}}
Suppose that $K$ is a satellite $L$-space knot. Then the JSJ graph of $K$ is  a rooted interval. 
\end{customthm}

\begin{customthm}{\ref{thm: jsj lo}}
Suppose that $K$ is a satellite knot which admits an irreducible rational surgery that is not $LO$. Then the JSJ graph of $K$ is a rooted interval. 
\end{customthm}

Theorem \ref{thm: jsj graph interval $L$-space knot} can also be deduced by combining the Baker-Motegi result that the pattern of an $L$-space satellite knot is braided in the pattern solid torus \cite{BMot19} and Lemma 3.1 of \cite{Be91}. 

We note that the hypothesis in Theorem \ref{thm: jsj lo} that the surgery is irreducible is necessary, though not particularly restrictive. A fundamental result of Scharlemann states that surgery on a satellite knot $K$ yields a reducible manifold if and only if $K$ is a cable knot and the surgery slope is the cabling slope (\cite{Sch90}; also see \S \ref{subsec: reducible toroidal surgery}). As such a surgery has a lens space summand, its fundamental group is not left-orderable. Hence Theorem \ref{thm: jsj lo} fails without the irreducibility assumption for a knot $K$ if and only if $K$ is a cable of a knot whose JSJ graph is not a rooted interval.

The Baker-Motegi result that the pattern of an $L$-space satellite knot is braided \cite{BMot19} implies that the pattern has winding number at least $2$. The following two theorems recover this latter fact using the techniques of the current paper (Corollary \ref{cor: w>1 nls}) and prove $LO$ analogues.  

\begin{customthm}{\ref{thm: w = 0}}
Let $K$ be a satellite knot with a pattern of winding number $0$. Let $r \in \mathbb Q$ and assume that $r$ is not the cabling slope if $K$ is a cable knot. Then $K(r)$ is $LO$ and $NLS$. 
\end{customthm}

\begin{customthm}{\ref{thm: w = 1}}
Let $K$ be a satellite knot with a pattern of winding number 1. Then for any $p\in \mathbb{Z}$, $K(p)$ is $LO$ and $NLS$. More generally, $K(p/(np \pm 1))$ is $LO$ and $NLS$ for any integers $p$ and $n$ such that $np \pm 1 \ne 0$. In particular, this holds for all rational surgeries $K(p/q)$ with $p = 1,2,3,4$, or $6$. 
\end{customthm}

As mentioned above, it is known that satellite knots of winding number $1$ are not $L$-space knots, and hence all rational surgeries on the knot are $NLS$. So by the $L$-space conjecture, it is expected that all rational surgeries on winding number $1$ satellite knots yield manifolds that are $LO$. Theorem \ref{thm: w = 1} shows that this holds for a set of slopes that is unbounded in both positive and negative directions. We refine this result in Theorem \ref{thm: farey graph} and use the refinement to produce examples of knots $K_k$, $k\in \mathbb{Z}^+$, for which $K_k(r)$ is $LO$ (and $NLS$) for any $r\in B_k(0)$, where $B_k(0)$ is the radius $k$ ball neighborhood of the $0$ slope in the Farey graph.

The key result that is used to prove Theorems \ref{thm: w = 0} and \ref{thm: split} is the fact that the meridional slope of a non-trivial knot $K$ is both $LO$-detected and $NLS$-detected \cite{BGH21} (see Theorem \ref{thm: meridional detn}). (We refer the reader to \S \ref{subsec: slope detection} for the definitions of a slope on $\partial X(K)$ being $LO$-, $NLS$- and $CTF$-detected.) The proof of Theorem \ref{thm: w = 1} uses the more general fact that any slope on the boundary of an integer homology solid torus at distance 1 from the longitude is $LO$- and $NLS$-detected. In \cite[Conjecture 1.6]{BGH21} we conjectured that the latter also holds for $CTF$-detection. In particular, 

\begin{conj}[cf. Conjecture 1.6 in \cite{BGH21}]
\label{conj: CTF meridional detection}
The meridional slope of a non-trivial knot is $CTF$-detected. 
\end{conj}

If Conjecture \ref{conj: CTF meridional detection} holds, then  the conclusions to Theorem \ref{thm: w = 0}, Theorem \ref{thm: split} and Corollary \ref{cor: nlo irred prime} hold with $NLS$ and $LO$ replaced by $CTF$. This provides a new approach to the problem of showing that any rational Dehn filling on a composite knot is $CTF$, as predicted by the $L$-space conjecture. Earlier work of Delman and Roberts \cite{DR21} studied this problem from the point of view of persistently foliar knots. (See Conjecture 1.6 of \cite{DR21}).  

In \cite{BGH21}, we showed that Conjecture \ref{conj: CTF meridional detection} holds for fibred  knots. This allows us to partially extend our results on $LO$ and $NLS$ surgeries to $CTF$ surgeries in the case that the companion knots are fibred. See Theorems \ref{thm: ctf winding zero fibered} and \ref{thm: ctf split fibered} for the precise statements; these are analogues of Theorems \ref{thm: w = 0} and \ref{thm: split} respectively. The following analogue of Theorem \ref{thm: jsj lo} is an immediate consequence of Theorem \ref{thm: ctf split fibered}.

\begin{customthm}{\ref{thm: ctf rooted interval}}
Suppose that $K$ is a fibered satellite knot which admits an irreducible rational surgery that is not $CTF$. Then the JSJ graph of $K$ is a rooted interval.
\end{customthm} 

Again, as in Theorem \ref{thm: jsj lo}, the only case where the irreducibility assumption is needed is when $K$ is a cable of a fibered knot whose JSJ graph is not a rooted interval.

Theorem \ref{thm: ctf rooted interval} implies the following result of Delman and Roberts, which is an immediate consequence of \cite[Theorem 6.1]{DR21}.

\begin{customcor}{\ref{cor: ctf composite fibered}}[Delman-Roberts]
All rational Dehn surgeries on a composite fibred knot are $CTF$. 
\end{customcor}

In addition to what we have mentioned above, it is known that if a satellite knot $K = P(K_0)$, with pattern $P$ and companion $K_0$, is an $L$-space knot, then both $K_0$ and $P(U)$ are $L$-space knots \cite{HRW1}. One therefore expects the $LO$ analogue of this to hold also. In this direction, we prove the following. 

\begin{customthm}{\ref{thm: all implies all}}
Let $K = P(K_0)$ be a satellite knot. Then for any $r \in \mathbb Q$, $K(r)$ is $LO$ if either
\begin{enumerate}[leftmargin=*] 
\setlength\itemsep{0.3em}
\item [{\rm (1)}] $K_0(s)$ is $LO$ for all $s \in \mathbb Q$ and $r$ is not the cabling slope of $K$ if $K$ is cabled, or
\item [{\rm (2)}] $P(U)(s)$ is $LO$ for all $s \in \mathbb Q$.
\end{enumerate}
\end{customthm}

Delman and Roberts call a knot $K$ {\em persistently foliar} if, for each rational slope $r$, there is a
co-oriented taut foliation meeting $\partial X(K)$ transversely in a foliation by curves of that slope.
Note that if a knot is persistently foliar, then $K(r)$ is $CTF$ for any $r\in \mathbb{Q}$, but the converse is unknown. 

\begin{customthm}{\ref{thm: satellite ctf}}
Suppose that $K$ is a satellite knot with a persistently foliar companion whose associated pattern has winding number $w\geq 1$. Then $K(r)$ is $CTF$ for any $r\in \mathbb{Q}$ which is not the cabling slope if $K$ is a cable knot. 
\end{customthm}

A special case of Theorem \ref{thm: satellite ctf} also recovers the following result of Delman and Roberts. 

\begin{customcor}{\ref{cor: dr composite}}[Delman-Roberts \cite{DR21}]
Each rational surgery on a composite knot with a persistently foliar summand is $CTF$. 
\end{customcor}

In \cite[Theorem 4.7]{RobertsSurfacebundle2}, Roberts shows that if $K$ is a fibred  knot whose monodromy has non-negative fractional Dehn twist coefficient, then for any slope $r\in (-\infty, 1)$, $K(r)$ is $CTF$. We combine her result with our methods to  show  the following:  

\begin{customcor}{\ref{cor: $L$-space ctf}}   
If $K$ is a positive satellite $L$-space knot with pattern of winding number $w$ then $K(r)$ is $CTF$ for each rational $r \in (-\infty, w^2]$. 
\end{customcor}

The fact that Conjecture \ref{conj: CTF meridional detection} is unproven is a major impediment to proving the $CTF$ analogues of our results in full generality. Question \ref{que: degree nonzero map} and the discussion following it explains another fundamental difficulty in proving the $CTF$ analogue of Theorem \ref{thm: all implies all}. For more results and discussion of $CTF$ Dehn surgery on satellite knots, see \S \ref{sec: surgery on satellite knots ctf}.

Finally, we present the following theorem regarding left-orderable $p/q$-surgeries when $p$ is small. 

\begin{thm}
\label{thm: small p}
If $K$ is a non-trivial knot, $p = 1$ or $2$, and $q \ne 0$, then the set rational numbers $p/q$ such $K(p/q)$ is not $LO$ is  contained in either $\{\pm1, \pm2 \}$ or $\{\varepsilon/2, 2\varepsilon/3, \varepsilon, 2\varepsilon \}$ for some $\varepsilon \in \{\pm1\}$.
\end{thm}

See Theorem \ref{thm: $LO$ surgery small p full} for a more detailed statement. We merely note here that although it is necessary to exclude $|p/q| = 1$ or $2$, the $L$-space Conjecture predicts that the potential exceptional slopes $\pm1/2$ and $\pm2/3$ do not arise. See Remark \ref{rem: mostly nec}.  

We now specialise to the case $p = 1$. As background, recall that Ozsv\'ath and Szab\'o conjectured that an irreducible $L$-space integer homology sphere is homeomorphic to either $S^3$ or the Poincar\'e homology  sphere $\Sigma(2, 3, 5)$. This is known to hold for integer homology $3$-spheres which arise as surgery on a knot $K$ by \cite{Ni07}, \cite{OS05-lens}, \cite{OS11}. Moreover, if $K$ is non-trivial and $K(r)$ is an $L$-space integer homology sphere, then up to replacing $K$ by its mirror image, $K$ is the right-handed trefoil knot $T(2, 3)$ and $r=1$. Hence, the $L$-space conjecture predicts that for $q\neq 0$, $K(1/q)$ is $LO$ unless $K = T(2, 3)$ and $q=1$ or $K = T(2, -3)$ and $q = -1$.  In this direction we have the following immediate corollary of Theorem \ref{thm: small p}.

\begin{cor}
    \label{cor: HF poincare LO}
At most two rational surgeries on a non-trivial knot yield integer homology $3$-spheres that are not $LO$, and the associated slopes are contained in either $\{-1, 1\}$, $\{-1, -1/2\}$, or $\{1/2, 1\}$. 
\end{cor}

\subsection*{Organisation of the paper} 
In Section \ref{sec: background} we review the various notions of slope detection, set our notational conventions for  satellite knots, and review some background results of Berge, Gabai and Scharlemann concerning surgery on satellite knots. In Section \ref{sec: $LO$ $NLS$ surgeries}, we use slope detection to prove results on $LO$ and $NLS$ surgeries on satellite knots. 
$CTF$ surgeries on satellite knots and related results are discussed in Section \ref{sec: surgery on satellite knots ctf}. In the final section \ref{sec: $LO$ surgery for p small}, we prove results on $LO$ surgeries on non-trivial knots for surgery coefficients $p/q$ with $p$ small. For hyperbolic knots (Section \ref{subsec: $LO$ surgery hyperbolic small p}), we use a fact about the Euler class of Fenley's asymptotic circle representation associated to a pseudo-Anosov flow, which is proved in the Appendix. 

\begin{acknowledgement*}
The authors wish to thank the referee for comments which led to an improvement of the exposition and Duncan McCoy and Patricia Sorya for helpful remarks related to Theorem \ref{thm: split}.
\end{acknowledgement*} 

\section{Preliminaries} 
\label{sec: background}
In this section we review various notions needed in the paper. 

\subsection{Slope detection and gluing results}
\label{subsec: slope detection}
A {\it slope} on a torus $T$ is an isotopy class of essential simple closed curves on $T$ or equivalently, a $\pm$-pair of primitive elements of $H_1(T)$. We will often simplify notation by denoting the slope corresponding to a primitive pair $\pm \alpha \in H_1(T)$ by $\alpha$.

The {\it distance} $\Delta(\alpha, \beta)$ between two slopes $\alpha, \beta$ on $T$ is the absolute value of the algebraic intersection number $|\alpha \cdot \beta|$ of their primitive representatives. Thus the distance is $0$ if and only if the slopes coincide and is $1$ if and only the representatives form a basis of $H_1(T)$. 

A {\it knot manifold} is a compact, connected, orientable, irreducible $3$-manifold $M$ with incompressible torus boundary. Given the  irreducibility of $M$, the condition that $M$ has incompressible torus boundary is equivalent to $M$ not being homeomorphic to $S^1 \times D^2$. 

The (rational) {\it longitude} of a knot manifold $M$ is the unique slope on $\partial M$ represented by a primitive class $\lambda_M \in H_1(\partial M)$ which is zero in $H_1(M; \mathbb Q)$. If $M$ is the exterior of a knot $K$ in a $3$-manifold $W$, there is a well-defined {\it meridional} slope represented by a primitive class $\mu \in H_1(\partial M)$ which is zero in $H_1(N(K))$. 

When $M$ is the exterior of a knot in the $3$-sphere, there are representatives $\mu, \lambda \in H_1(\partial M)$ of the meridional and longitudinal slopes which are well-defined up to simultaneous sign change. This yields a canonical identification of the set of slopes on $\partial M$ and $\mathbb Q \cup \{\frac{1}{0}\}$ via $\pm(p \mu + q \lambda) \longleftrightarrow p/q$, where $p$ and $q$ are coprime integers.   

We say that a slope $\alpha \in H_1(\partial M)$ on the boundary of a knot manifold $M$ is: 

\begin{itemize} 
\setlength\itemsep{0.3em}
\item {\it $CTF$-detected} (in $M$) if there is a co-oriented taut foliation $\mathcal{F}$ on $M$ which intersects $\partial M$ transversely and is such that $\mathcal{F}\cap \partial M$ is a foliation without Reeb annuli which contains a closed leaf of slope $\alpha$ \footnote{This is equivalent to the definition given in \cite[Definition 5.1]{BGH21}. See \cite[Proposition 4.3.2]{HH81} for the relevant background on foliated tori.}. 

\item {\it $LO$-detected} (in $M$) if $N \cup_f M$ has a left-orderable fundamental group, where $N$ is the twisted $I$-bundle over the Klein bottle and $f: \partial N \xrightarrow{\; \cong \;} \partial M$ identifies the rational longitude of $N$ with $\alpha$ \footnote{The definition for $LO$-detection here is in fact a characterization of the $LO$-detected slopes shown in \cite[Theorem 1.4]{BC2} (also see \cite[Corollary 6.12]{BGH21}). The original definition of $LO$-detection is given in \cite[Definition 7.5]{BC2} and \cite[Definition 6.2]{BGH21}.}.

\item {\it $NLS$-detected} (in $M$) if it lies in the closure of the set of slopes whose associated Dehn filling is not an $L$-space (\cite[\S 7.2]{RR17}). Here the closure is taken inside the set of all rational slopes, which is homeomorphic to $\mathbb{Q}P^1$.
\end{itemize}

\begin{remark}[Slope detection and Dehn filling]
\label{rem: detection and Dehn fillings}
It follows by definition that if the $\alpha$-Dehn filling of $M$ is $NLS$, then $\alpha$ is $NLS$-detected. Similarly, if $M(\alpha)$ is $LO$ then $\alpha$ is $LO$-detected \cite[Corollary 8.3]{BC2}. The converse statements are false in general and therefore the condition that a slope $\alpha$ be $NLS$- or $LO$-detected is {\it strictly weaker} than requiring that the Dehn-filled manifold $M(\alpha)$ be $LO$ or $NLS$. For example, even though the meridional slope $\mu$ of a non-trivial knot is both $LO$-detected and $NLS$-detected in the knot's exterior $M$ (cf. Theorem \ref{thm: meridional detn}), the $\mu$-Dehn filling of $M$ is $S^3$ and therefore neither $LO$ nor $NLS$. 

It is not known whether a Dehn filling $M(\alpha)$ of a knot manifold $M$ being $CTF$ implies that the slope $\alpha$ is $CTF$-detected, though this is thought to be true (and is true if $M$ is a Seifert fibre space \cite{BC17}). If so, the condition that $\alpha$ is $CTF$-detected will be strictly weaker than the condition that $M(\alpha)$ be $CTF$: the meridional slope $\mu$ of a non-trivial fibred knot is $CTF$-detected (cf. Theorem \ref{thm: meridional detn}) though the $\mu$-Dehn filling is $S^3$ and therefore not $CTF$. 
\end{remark}

\begin{prop}
\label{prop: longitude detected}
The longitudinal slope of a knot manifold is $\ast$-detected, where $*$ denotes either $NLS$, $CTF$, or $LO$. 
\end{prop} 

\begin{proof} The proposition is obvious for $NLS$-detection, since $L$-spaces have finite first homology. That it is $CTF$-detected follows from the main result of \cite{Gab83} and $LO$-detected from \cite[Example 6.3]{BGH21}.   
\end{proof}

Slope detection is well-adapted to understanding when a union $W = M_1 \cup_\partial  M_2$ of two knot manifolds has property $\ast$, where ``$\ast$" stands for either $NLS$,  $LO$, or $CTF$. More precisely, combining work of Hanselman, Rasmussen and Watson \cite[Theorem 13]{HRW1}, Boyer and Clay \cite[Theorem 1.3]{BC2}, and  \cite[Theorem 5.2]{BGH21} yields the following gluing theorem. 

\begin{thm}[Theorem 1.1 in \cite{BGH21}]
\label{thm: * gluing}
Suppose that $W = M_1 \cup_f M_2$ where $M_1, M_2$ are knot manifolds and $f: \partial M_1 \xrightarrow{\; \cong \; } \partial M_2$. 
If $f$ identifies a $\ast$-detected slope on the boundary of $M_1$ with a $\ast$-detected slope on the boundary of $M_2$, then $W$ has property $\ast$. 
\end{thm}

\begin{thm}[Theorem 1.3 in \cite{BGH21}]
\label{thm: meridional detn} 
Let $M$ be an integer homology solid torus knot manifold. Then each slope of distance $1$ from $\lambda_M$ is $LO$-detected and $NLS$-detected. If $M$ fibres over the circle, any slope of distance $1$ from $\lambda_M$ is also $CTF$-detected. 
\end{thm}

One can similarly define $\ast$-detection for multislopes in the case of manifolds with multiple toral boundary components. See, respectively, \S 4.3, \S 5.3, \S 6.6 of \cite{BGH21} for the definitions of $NLS$, $CTF$, and $LO$ multislope detection. There is an associated gluing theorem (\cite[Theorem 7.6]{BGH21}): If $\ast$-detected multislopes coincide along matched boundary components, then the resulting manifold has the property $\ast$. In this article we use multislope detection and gluing to prove Theorem \ref{thm: split}.

\subsection{Conventions on satellite knots}
\label{subsec: conventions}
Throughout this paper, we will express a satellite knot $K$ with pattern $P$ and companion $K_0$ as $P(K_0)$. In this setting, we can write 
$$S^3 = V \cup_T X_0$$ 
where $V$ is a solid torus containing $P$, $X_0$ is the exterior of $K_0$, and $T = \partial V = \partial X_0$. We require that $K_0$ be a non-trivial knot and that $P$ is neither isotopic to the core of $V$ nor contained in a $3$-ball in $V$.  As such, $T$ is essential in $X(K)$. 

  Conversely, any essential torus $T_1$ in $X(K)$ bounds the exterior of a non-trivial knot $X(K_1) \subset X(K)$ to one side and a solid torus $V_1 \subset S^3$ containing $K$ to the other. Thinking of $K$ as a pattern knot $P_1 \subset V_1$ allows us to express $K$ as a satellite knot $P_1(K_1)$. Since it is possible for $X(K)$ to contain non-isotopic essential tori, it is possible for $K$ to be expressed as a satellite knot in different ways. 

If $P$ and $P_0$ are patterns we define the {\it satellite pattern} $P(P_0)$ in the obvious way. 

When we express a satellite knot as $K = P(K_0)$, the {\em winding number} of the pattern $P$ in $V$ is the non-negative integer $w = w(P)$ for which 
$$\text{image}(H_1(P) \to H_1(V)) = wH_1(V)$$ 
This coincides with the {\em winding number} of $T = \partial X_0$ in the exterior $X(K)$ of $K$, which is defined to be the non-negative integer $w = w(T)$ for which
\begin{equation}
    \label{equ: winding number T}
    \text{image}(H_1(T) \to H_1(X(K))) = wH_1(X(K))
\end{equation}
In what follows, we will express our results in terms of the winding numbers of either patterns or tori depending on our focus.

Given an expression $P(K_0)$ of $K$ as a satellite knot, let $N(P)$ denote a regular neighborhood of $P$ in $V$. We will use $M$ to denote the exterior
of $P$ in $V$: 
\begin{equation} 
\label{eqn: pattern exterior} 
M = V\setminus \text{int}(N(P)) = X(K) \setminus \text{int}(X_0)
\end{equation}

Note that $M$ is irreducible, has incompressible boundary, and is not a product $T^2 \times I$. 

Let $\mu_0, \lambda_0 \in H_1(\partial X_0)$ denote meridional and longitudinal classes of $K_0$ and $\mu, \lambda \in H_1(\partial X(K))$ denote those of $K$. 
These slopes can be oriented so that in $H_1(M)$ we have 
\begin{equation}
    \label{equ: meridians}
    \mu_0 = w \mu \quad \text{and} \quad w \lambda_0 = \lambda
\end{equation}

For a coprime pair $p$ and $q$ with $p \geq 0$, define $M(p/q)$ to be the $p \mu + q \lambda$ Dehn filling of $M$ along $\partial X(K)$. Equivalently, $M(p/q)$ is the $p/q$-Dehn surgery on the knot $P$ in $V$. We often denote $M(p/q)$ by $P(p/q)$. Similarly, we use $K(p/q)$ to denote the $p \mu + q \lambda$ Dehn filling of $X(K)$ along $\partial X(K)$. Then 
\begin{equation}
\label{eqn: surgery decomp} 
K(p/q) = P(p/q) \cup_T X_0 
\end{equation}
From (\ref{equ: meridians}), it follows that $p \mu_0 + qw^2 \lambda_0 = w (p\mu + q\lambda)$ is null-homologous in $P(p/q)$. Hence the longitudinal slope $\lambda_{P(p/q)}$ of $P(p/q)$ is given by 
\begin{equation}
\label{eqn: lambda} 
\lambda_{P(p/q)} = \frac{1}{\gcd(p, w^2)} \left(p \mu_0 + qw^2 \lambda_0 \right) 
\end{equation}
It follows that the element of $\mathbb Q\cup \{\frac{1}{0}\}$ corresponding to 
$\lambda_{P(p/q)}$ as a slope on $\partial X_0$ is 
\begin{equation}
\label{eqn: slope as a rational number}
\left\{ 
\begin{array}{cl}
1/0 & \mbox{ if } w = 0 \\
p/qw^2 & \mbox{ if } w \ne 0
\end{array} \right. 
\end{equation}

\begin{lemma} 
\label{lem: w=0 longitude}
Suppose $r \in \mathbb Q$. If $w=0$, then the longitudinal slope $\lambda_{P(r)}$ of $P(r)$ is the meridional slope of $X_0$, and therefore is $LO$-detected and $NLS$-detected in $X_0$. If $X_0$ is fibered, then $\lambda_{P(r)}$ is also $CTF$-detected in $X_0$.
\end{lemma}
\begin{proof}
The first claim follows from Equation (\ref{eqn: slope as a rational number}). The second and third then follow from Theorem \ref{thm: meridional detn}.
\end{proof}

\begin{lemma}
\label{lem: r/w^2 detected}
Let $K = P(K_0)$ be a satellite knot with pattern $P$ and companion $K_0$. Fix $r \in \mathbb Q$. If $T = \partial X_0$ is incompressible in $K(r)$ and $\lambda_{P(r)}$ is $*$-detected  in $X_0$, then $K(r)$ has property $*$, where $*$ is either $LO, NLS$, or $CTF$.
\end{lemma}
\begin{proof}
By (\ref{eqn: surgery decomp}), we have $K(r) = P(r) \cup_T X_0$. Since $T$ is incompressible in $K(r)$ by assumption, both $P(r)$ and $X_0$ are knot manifolds. By Proposition \ref{prop: longitude detected}, the longitudinal slope $\lambda_{P(r)}$ of $P(r)$ is $\ast$-detected in $P(r)$. As  we have assumed that $\lambda_{P(r)}$ is $\ast$-detected in $X_0$, Theorem \ref{thm: * gluing} implies that $K(r)$ has property $\ast$.
\end{proof}

\subsection{Conventions on cable knots} 
Recall that given coprime integers $m$ and $n$ with $|m|\geq 2$ (so $n \ne 0$), the $(m,n)$-cable $C_{m,n}(K_0)$ of a non-trivial knot $K_0$ is the satellite knot with companion $K_0$ and pattern $C_{m,n}$, the $(m, n)$ torus knot $T(m,n)$ standardly embedded in $V$ with winding number $|m|$. Since reversing the orientation of $C_{m,n}(K_0)$ changes the signs of $m$ and $n$ simultaneously, without loss of generality we can assume below that $m\geq 2$ and $n\neq 0$. The slope $mn$ is called the {\it cabling slope} of the cable knot $C_{m,n}(K_0)$. The exterior  of the pattern knot $C_{m,n}$ in $V$ is an {\it $(m, n)$-cable space} and is a Seifert fibre space whose base orbifold is an annulus with a single cone point, whose order is $m$.

If $P_0$ is a pattern we can define the satellite pattern $C_{m,n}(P_0)$. We say a pattern is {\it cabled} if it is either $C_{m,n}$ or $C_{m,n}(P_0)$ or some pattern $P_0$. 

Though there are potentially many ways to express a knot as a satellite knot, the following lemma shows that there is only one way of expressing it as a cable knot. In other words, if $C_{m,n}(K_0) = C_{m',n'}(K_1)$, then $m' = m, n' = n$ and $K_0 = K_1$. Hence there is a unique cabling slope on $\partial X(K)$ for a cable knot, and we refer to this slope as {\it the cabling slope} of $K$.

\begin{lemma}
\label{lemma: cabled implies never 1}
Suppose that $K = C_{m,n}(K_0)$. Then any essential torus in $X(K)$ is isotopic into the exterior $X_0$ of $K_0$. Consequently,
\begin{enumerate}[leftmargin=*]
\setlength\itemsep{0.3em}
\item[{\rm (1)}] If $K = P(K_1)$ then $P$ is cabled. More precisely, if $K = P(K_1)$, then either
\begin{itemize}
\item[{\rm (a)}] $P = C_{m,n}$ and $K_0 = K_1$, or
\item[{\rm (b)}] $P = C_{m,n}(P_1)$ for some pattern $P_1$ and $K_0 = P_1(K_1)$;
\end{itemize}
\item[{\rm (2)}] $K = C_{m,n}(K_0)$ is the unique realisation of K as a cable knot and therefore $mn$ is the unique slope on $\partial X(K)$ which is the slope of a cabling of $K$; 
\item[{\rm (3)}] the winding number of any essential torus $T$ in $X(K)$ is divisible by $m$. 
\end{enumerate}
\end{lemma}

\begin{proof}
Write $X(K) = M \cup_{T_0} X_0$, where $M$ is the exterior of $C_{m,n}$ in the solid torus $V$,  so $M$ is a cable space and $\partial M = \partial X(K) \cup T_0$. Since the only Seifert pieces of the JSJ decomposition of the exterior of a non-trivial knot are cable spaces, composing spaces, or torus knot exteriors \cite[Lemma VI.3.4]{JS79}, $M$ is a piece of the JSJ decomposition of $X(K)$.

Let $T$ be an essential torus in $X(K)$. By the characteristic property of the JSJ decomposition, $T$ can be isotoped into a Seifert fibred piece of the JSJ decomposition of $X(K)$. Since the cable space $M$ contains no essential tori, it follows that $T$ can be istoped into $X_0$.

Suppose $K = P(K_1)$. This gives a decomposition $X(K) = M_1 \cup_{T_1} X_1$, where $M_1$ is the exterior of $P$ in the solid torus $V_1$ bounded by $T_1$ in $S^3$, and $X_1$ is the exterior of $K_1$. By the above, $T_1$ can be isotoped into $X_0$. If $T_1$ is parallel to $T_0$ we have conclusion (a). Otherwise, we have $V \subset V_1$, where the core of $V$ is a pattern $P_1$ in $V_1$. Hence $P = C_{m,n}(P_1)$. Also, $X_0 = M_1 \cup_{T_1} X_1$, where $M_1$ is the exterior of $P_1$ in $V_1$, so $K_0 = P_1(K_1)$, giving (b).

(2) follows by taking $P = C_{m',n'}$ and noting that here conclusion (b) is impossible since cable spaces are atoroidal.

For (3), isotope $T$ into $X_0$ so that $T\subset X_0 \subset X(K)$.  Since the winding number of $T_0 = \partial X_0$ in $X(K)$ is $m$, by (\ref{equ: winding number T}) we have $${\rm image}(H_1(X_0)\rightarrow H_1(X(K))) = {\rm image}(H_1(T_0)\rightarrow H_1(X(K))) = m H_1(X(K))$$ 
Then as the homomorphism $H_1(T) \to H_1(X(K))$ factors through $H_1(X_0)$, $m$ must divide the winding number of $T$ in $X(K)$. 
\end{proof}

\begin{remark}
Suppose the JSJ graph of a satellite knot $K$ is a rooted interval; equivalently, all the pieces of the JSJ decomposition of $X(K)$ have exactly two boundary compoments except the innermost piece, which has a single boundary component. Each of the former is therefore either a cable space or hyperbolic, and the latter is the exterior of a knot $K_0$ which is either a torus knot or hyperbolic. Suppose the JSJ graph has $k \ge 1$ edges, with corresponding JSJ tori $T_i$, $1 \le i \le k$. Let $P_k(...(P_1(K_0))...)$ be the associated iterated satellite expression of $K$. As argued in the proof above, every essential torus in $X(K)$ is isotopic to some $T_i$, which implies that the given length $k$ iterated satellite expression of $K$ is unique.
\end{remark}

\subsection{Reducible and toroidal surgery on satellite knots}
\label{subsec: reducible toroidal surgery}
Let $K$ be a satellite knot with pattern $P$, a knot in the solid torus $V$, and companion knot $K_0$. We adhere to the notations set in \S \ref{subsec: conventions}.  

The manifold $K(p/q)$ has the following decomposition given in (\ref{eqn: surgery decomp}):
\begin{displaymath}
K(p/q) = P(p/q) \cup_T X_0, 
\end{displaymath}
where $P(p/q)$ is the $p/q$-Dehn surgery on $P$ in $V$ and $X_0$ is the exterior of $K_0$. Many of the results in this article are obtained by applying Theorem \ref{thm: * gluing} to this decomposition. This requires that $P(p/q)$ be irreducible and $T = \partial P(p/q)$ to be incompressible in $P(p/q)$\footnote{$T$ is incompressible in $X_0$ by assumption.}. We remark that if $P(p/q)$ is irreducible but has a compressible boundary, then $P(p/q)$ is a solid torus.

\begin{thm}{\rm (Gabai \cite[Theorem 1.1]{Gab89}, Scharlemann \cite{Sch90})}
\label{thm: sch red}
Suppose that $P$ is a pattern in a solid torus $V$. If $r \in \mathbb Q$, then exactly one of the following three possibilities arises.
\begin{enumerate}[leftmargin=*]
\setlength\itemsep{0.3em}
\item [{\rm (a)}] $P(r)$ is a solid torus and $P$ is either a $0$-bridge braid or a $1$-bridge braid in $V$.
\item [{\rm (b)}]$P(r)$ is reducible, $P$ is cabled in $V$, and $r$ is the slope of the cabling annulus.
\item [{\rm (c)}]$P(r)$ is irreducible with incompressible boundary. 
\end{enumerate}
\end{thm}

In the statement of Theorem \ref{thm: sch red}(a), a {\em $0$-bridge braid} in a solid torus $V$ is a torus knot standardly embedded in $V$, in other words, a cable pattern $C_{m,n}$. We refer the readers to \cite[\S 2]{Gab90} for details on $1$-bridge braids in $V$, though we note that they are parameterised up to orientation-preserving ambient homeomorphism by triples $(w, b, t)$, where $w \geq 3$ is the number of the strands of the braid, which equals its winding number in $V$, and $1 \leq b, t \leq w-2$.

\begin{lemma}
\label{lemma: 1 bridge winding number}
Suppose that $P$ is a $1$-bridge braid in a solid torus $V$ and there exists a non-trivial Dehn-surgery on $P$ that produces $S^1\times D^2$. Then $P$ has winding number $w \geq 4$ and if $w = 4$, $P$ is a cable on a $0$-bridge braid. 
\end{lemma}

\begin{proof}
If $w = 3$, then $b = 1$ since $1\leq b \leq w -2$. In this case Gabai showed that $P$ is a $(2, n)$-cable on a $0$-bridge braid in $V$ for some odd $n$ (\cite[Example 3.7]{Gab90}). But then $2$ divides $w$, which is a contradiction. Hence $w \geq 4$. 

If $w = 4$, then $b\leq 2$. 

Suppose first that $b = 2$. Then using the notation of \cite[Example 3.8]{Gab90}, $P$ has Type $1, 2$ or $3$. Type 1 is ruled out since in this case $w \equiv 0$ (mod $3$). Type 3 is ruled out since in this case there are integers $r, s \geq 1$ such  that $w = 3r + 2s \geq 5$. Thus $P$ has Type $2$ and in this case there are integers $r, s \geq 1$ such $4 = w = 3r + s$. Hence $r = s = 1$. But then Proposition 3.9 of \cite{Gab90} shows that no non-trivial Dehn-surgery on $P$ produces $S^1\times D^2$, which contradicts our hypotheses. 

Finally, if $w = 4$ and $b = 1$,  Example 3.7 of \cite{Gab90} shows that $P$ is a cable on a $0$-bridge braid, which completes the proof. 
\end{proof}

Here is an immediate application of Theorem \ref{thm: sch red}. 
\begin{cor}[Corollary 4.5 in \cite{Sch90}] 
\label{cor: sch red sat}
If $K$ is a satellite knot and $K(r)$ is reducible for some rational number $r$, then $K$ is a cable knot and $r$ is the cabling slope. 
\end{cor}

\begin{cor}
\label{cor: 1 implies irreducible} 
If $K = P(K_0)$ is a satellite knot where $P$ is a winding number $1$ pattern, then $K$ is not a cable knot and hence all surgeries on $K$ are irreducible. 
\end{cor}
\begin{proof}
Since the winding number of $K$ is $1$, Lemma \ref{lemma: cabled implies never 1}(3) shows that $K$ cannot be a cable knot. Hence all surgeries on $K$ are irreducible by Corollary \ref{cor: sch red sat}. 
\end{proof}

Theorem \ref{thm: C(T)} below is a combination of well-known results from Berge, Gabai, Gordon, and Scharlemann. Assume that $K = P(K_0)$ is a satellite knot with $X(K) = M \cup_T X_0$ as above. 
Set 
$$\mathcal{C}(T) = \{ r \in \mathbb{Q} \cup \left\{1/0\right\} \; | \; T \mbox{ compresses in } K(r)\}$$ 
Since $K(1/0) = S^3$, we have $\frac{1}{0} \in \mathcal{C}(T)$.

\begin{thm}{\rm (Berge \cite{Be91}, Gabai \cite{Gab89, Gab90}, Gordon \cite{Go83} and Scharlemann \cite{Sch90})}
\label{thm: C(T)}
Let $K = P(K_0)$ be a satellite knot with pattern $P$ of winding number $w$ and companion $K_0$. Set $T = \partial X_0$ and suppose that there is a non-meridional slope in $\mathcal{C}(T)$, i.e., $\mathcal{C}(T) \neq \{1/0\}$. Then$:$ 
\begin{enumerate}[leftmargin=*]
\setlength\itemsep{0.3em}
\item[{\rm (1)}] If $K$ is not a cable knot, then $P$ is a $1$-bridge braid with winding number $w \geq 5$. Further, 
there is an integer $a$ such that 
$$\mathcal{C}(T) \subseteq \{1/0, a, a+1\} $$
and $K(r) = K_0(r/w^2)$ for $r\in \mathcal{C}(T)$
\item[{\rm (2)}]  If $K = C_{m, n}(K_0)$ for coprime integers $m \geq 2$ and $n$,  then 
\medskip 
\begin{enumerate}
\setlength\itemsep{0.3em}
    \item[{\rm (a)}]  $\mathcal{C}(T)$ is the set of slopes of distance at most $1$ from the cabling slope $mn$. %
    \item[{\rm (b)}] If $r = mn \in \mathcal{C}(T)$, then $K(r) \cong K_0(n/m) \# L(m, n)$.
    \item[{\rm (c)}] For $r = p/q \in \mathcal{C}(T)$ that is distance equal to $1$ from the cabling slope $($that is, $p/q = mn \pm \frac{1}{b}$ for some $b\in \mathbb{Z})$, we have  $K(p/q) \cong K_0(p/m^2q)$.
\end{enumerate}
\end{enumerate}
\end{thm}

\begin{proof}
Suppose that $K = P(K_0)$ is not a cable knot and $T$ compresses in $K(r)$ for some $r \in \bQ$. As in (\ref{eqn: surgery decomp}), we write $K(r) = P(r) \cup_T X_0$. Since $T$ is incompressible in $X_0$, it must compress in $P(r)$. Then by Theorem \ref{thm: sch red}, we have that $P$ is a $1$-bridge braid and $P(r) \cong S^1 \times D^2$. Lemma \ref{lemma: 1 bridge winding number} then shows that $w\geq 5$. Further, \cite[Lemma 3.2]{Gab90} shows that $\mathcal{C}(T)$ contains at most two consecutive integer slopes in addition to the meridional slope $\frac{1}{0}$. Finally, since $P(r)$ is a solid torus, $K(r) = P(r)\cup X_0$ is a Dehn filling on $X_0$, and the slope of this surgery is $r/w^2$ by (\ref{eqn: slope as a rational number}). This proves (1). 

For (2), see Section 7 of \cite{Go83}.
\end{proof}

\begin{lemma}
\label{lemma: incomp torus}
Let $K = P(K_0)$ where $P$ has winding number $w$ and $X_0$ is the exterior of $K_0$. Then $K(w^2)$ is irreducible and $\partial X_0$ is incompressible in $K(w^2)$.     
\end{lemma}

\begin{proof}
The first conclusion follows from the second by Corollary \ref{cor: sch red sat} and Theorem \ref{thm: C(T)}. Also by that theorem, if $\partial X_0$ compresses in $K(r)$ and $r$ is integral then either $P = C_{m,n}$ and $r = mn$ or $mn \pm 1$, or $P$ is a 1-bridge braid and $P(r)$ is a solid torus. 

In the first case $r \ne m^2$. In the second case, by \cite{Gab90} $r = a$ or $a + 1$ where $a = tw + b$ and $1 \le b,t \le w-2$. Then $a+1 \le w(w-2) + (w-2) + 1 = w^2 - w - 1 < w^2$.
\end{proof}

\section{Left-orderability, non-\texorpdfstring{$L$}{L}-spaces and surgeries on satellite knots} 
\label{sec: $LO$ $NLS$ surgeries}

\subsection{\texorpdfstring{$LO$}{LO} and \texorpdfstring{$NLS$}{NLS} surgery on satellite knots and winding numbers}
\label{subsec: $LO$ winding numbers}

Baker and Motegi showed that the pattern of an $L$-space satellite knot must be braided \cite{BMot19}, so its winding number is at least $2$. In this section, we reprove this fact using slope detection and gluing techniques. See Corollary \ref{cor: w>1 nls} below. An identical argument shows that the winding number of any pattern of a satellite knot which admits irreducible non-$LO$ Dehn fillings must be strictly positive.  

\begin{thm}
\label{thm: w = 0} 
Let $K$ be a satellite knot with a pattern of winding number $0$. Let $r \in \mathbb Q$ and assume that $r$ is not the cabling slope if $K$ is a cable knot. Then $K(r)$ is $LO$ and $NLS$.  
\end{thm}

\begin{proof}
Let $K = P(K_0)$ where the winding number of $P$ is zero, and $K(r) = P(r) \cup_{T} X_0$ as before. Since $w(P) = 0$, $T$ is essential in $K(r)$ by Gabai \cite[Corollary 2.5]{Gab2}. Since $r$ is not the cabling slope, by Theorem \ref{thm: sch red}, $P(r)$ is irreducible. Hence by Proposition \ref{prop: longitude detected}, the longitudinal slope of $P(r)$ is $LO$- and $NLS$-detected in $P(r)$. By Lemma \ref{lem: w=0 longitude}, the longitudinal slope of $P(r)$ is also $LO$- and $NLS$-detected in $X_0$. Then the result follows from Lemma \ref{lem: r/w^2 detected}. 
\end{proof}

\begin{remark}
\label{rem: cabling slope assumption}
It is possible for a cable knot to admit an essential winding number zero torus in its exterior. This is the case, for instance, when $K$ is a cable of a Whitehead double of a knot: $K$ can be expressed as a cable satellite knot and also as a satellite knot with a pattern of winding number $0$. Thus it is necessary to exclude the possibility that $r$ is a cabling slope in the statement of Theorem \ref{thm: w = 0}. 
\end{remark}

Next we consider the case of a satellite knot $K$ which has a pattern of winding number $1$. 

\begin{thm}
\label{thm: w = 1}
Let $K$ be a satellite knot with a pattern of winding number $1$. Then for any $p \in \mathbb{Z}$, $K(p)$ is $LO$ and $NLS$. More generally, $K(p/(np \pm 1))$ is $LO$ and $NLS$ for any integers $p$ and $n$ such that $np \pm 1 \ne 0$. In particular, this holds for all rational surgeries $K(p/q)$ with $p = 1,2,3,4$, or $6$. 
\end{thm}

\begin{proof}
Let $K=P(K_0)$ where the winding number $w$ of the pattern $P$ is $1$.
Write $X(K) = M \cup_{T} X_0$ where $X_0$ is the exterior of $K_0$, $T = \partial X_0$ and $\partial M = \partial X(K) \cup T$ as before. Then $K(p/q)$ can be decomposed as $P(p/q) \cup_T X_0$.

Since $w= 1$, $K(p/q)$ is irreducible for all $p/q \in \mathbb Q$ by Corollary \ref{cor: 1 implies irreducible}. Further, by Theorem \ref{thm: C(T)}(1), $T$ remains incompressible in $K(p/q)$. 

Another consequence of our assumption that $w = 1$ is that $P(p/q)$ is an integer homology solid torus for each $p/q$. Then by Theorem \ref{thm: meridional detn}, any slope on $T$ that is distance $1$ from the longitudinal slope $\lambda_{P(p/q)}$ of $P(p/q)$ is $LO$- and $NLS$-detected. Since $X_0$ is also an integer homology solid torus, any slope on $T$ that is distance $1$ from its longitudinal slope $\lambda_0$ is also $LO$- and $NLS$-detected. As $T$ remains incompressible in $K(p/q)$, by Theorem \ref{thm: * gluing}, we have that $K(p/q)$ is $LO$ and $NLS$ as long as there is a slope $\alpha$ on $T$ satisfying that $\Delta(\alpha, \lambda_0) = \Delta(\alpha, \lambda_{P(p/q)}) = 1$. 

But $\Delta(\alpha, \lambda_0) = 1$ if and only if $\alpha$ corresponds to $\mu_0 + n \lambda_0$ for some $n \in \mathbb Z$, where $\mu_0$ is the meridian of $K_0$. On the other hand, $\lambda_{P(p/q)} = p \mu_0 + q \lambda_0$ by (\ref{eqn: slope as a rational number}). So $1 = \Delta(\mu_0 + n \lambda_0, \lambda_{P(p/q)}) = |np - q|$ if and only if $q = np \pm 1$. 

In the case that $p\in \pm \{1, 2, 3, 4, 6\}$, each integer coprime with $p$ is congruent to $\pm 1$ (mod $p$), which yields the final claim of the theorem.
\end{proof}

The set of $L$-space surgery slopes of an $L$-space knot of genus $g$ is the set of slopes in either $[2g-1, \infty]$ or $[-\infty, -(2g-1)]$ \cite{OS05-lens,OS11} and therefore the fact that the winding numbers of $L$-space satellite knots must be strictly bigger than $1$ follows immediately from Theorems \ref{thm: w = 0} and \ref{thm: w = 1}.

\begin{cor}
\label{cor: w>1 nls} 
If $K$ is a satellite $L$-space knot then any essential torus in the exterior of $K$ has winding number $2$ or more. 
\end{cor}

The $L$-space Conjecture predicts that all rational surgeries on a satellite knot with winding number $1$ are $LO$. Theorem \ref{thm: w = 1} shows that this holds for a set of slopes that is unbounded in both positive and negative directions, but our present lack of knowledge about the set of $LO$-detected slopes on the boundary of an integer homology solid torus prevents us from proving this for all slopes. However, by considering iterated winding number $1$ satellites we at least get knots $K$ with successively larger sets of slopes $r$ for which we can show that $K(r)$ is $LO$.  

Let $\bQ_{\infty} = \mathbb{Q}\cup \{1/0\}$. To state the precise result, recall that the Farey graph is the graph with vertices $\bQ_{\infty}$ and an edge between $r,r' \in \bQ_{\infty}$ if and only if $\Delta(r,r') = 1$. We denote by $d_{FG}(r,s)$ the distance in the Farey graph between $r,s \in \bQ_{\infty}$. 

\begin{thm}
    \label{thm: farey graph}
  Let $K$ be an iterated satellite knot $P_{k-1}(...(P_1(K_0))...)$, where $k \ge 2$ and the winding number $w(P_i) = 1$ for all $i$. Then for $r \in \bQ$, $K(r)$ is $LO$ and $NLS$ if $d_{FG}(0,r) \le k$.  
\end{thm}

\begin{proof}
We prove the statement for $LO$. For $NLS$, the argument is identical; one simply replaces $LO$ by $NLS$ in the proof below.

For $k \ge 0$ let $B_k = \{r \in \bQ_{\infty} : d_{FG}(0,r) \le k \}$. Thus $B_0 = \{0\}$, $B_1 = \{1/n : n \in \bZ\}\cup \{0\}$, and $B_2 = \{p/(np\pm1) : p,n \in \bZ\}$.

Let $K$ be as in the theorem. We must show that $K(r)$ is $LO$ for all $r \in B_k \setminus \{1/0\}$. We proceed by induction on $k$; the statement holds for $k = 2$ by Theorem \ref{thm: w = 1}.

So suppose that $k \ge 3$, and that the claim holds for $k-1$. To simplify the notation, write $P = P_{k-1}$ and $K' = P_{k-2}(...(P_1(K_0))...)$, so $K = P(K')$. The inductive hypothesis says that $K'(r')$ is $LO$ for all $r' \in B_{k-1} \setminus \{1/0\}$, and hence all slopes in $B_{k-1} \setminus \{1/0\}$ are $LO$-detected in $X(K')$ by Remark \ref{rem: detection and Dehn fillings}. Also, $1/0$ is $LO$-detected in $X(K')$ by Theorem \ref{thm: meridional detn}. Therefore any $r' \in B_{k-1}$ is $LO$-detected in $X(K')$.

Let $r \in B_k\setminus \{1/0\}$ represent a slope on $\partial X(K)$. Then $K(r) = P(r) \cup_{T} X(K')$, where $P(r)$ is irreducible and $T = \partial P(r) = \partial X(K')$ is incompressible in $K(r)$ by Corollary \ref{cor: 1 implies irreducible} and Theorem \ref{thm: C(T)}. To prove that $K(r)$ is $LO$, we will show that there exists a slope $r' \in B_{k-1}$ that is $LO$-detected in $P(r)$. Since $r'\in B_{k-1}$, it is also $LO$-detected in $X(K')$ by our discussion above. Then the conclusion follows from Theorem \ref{thm: * gluing}.

Since $w(P) = 1$, $P(r)$ is an integer homology solid torus with longitudinal slope $r$ on $\partial X(K')$ by (\ref{eqn: slope as a rational number}). Since $r \in B_k\setminus \{1/0\}$, there exists $r' \in B_{k-1}$ such that $\Delta(r,r') \le 1$. The slope $r'$ is $LO$-detected in $P(r)$ by Theorem \ref{thm: meridional detn} and Proposition \ref{prop: longitude detected}.  This completes the proof.
\end{proof}

\subsection{\texorpdfstring{$LO$}{LO} and \texorpdfstring{$NLS$}{NLS} surgery on satellite knots from patterns and companions}
\label{subsec: $LO$ from pattern and companion}
It is shown in \cite{HRW1} that if $K = P(K_0)$ is a satellite $L$-space knot, then both $K_0$ and $P(U)$ are $L$-space knots, where $U$ is the unknot. Equivalently, if either $K_0$ or $P(U)$ is not an $L$-space knot, then $K = P(K_0)$ is not an $L$-space knot. In this section, we prove the following analogous result for left-orderability. 

\begin{thm} 
\label{thm: all implies all} 
Let $K = P(K_0)$ be a satellite knot. Then for any $r \in \mathbb Q$, $K(r)$ is $LO$ if either
\begin{enumerate}[leftmargin=*] 
\setlength\itemsep{0.3em}
\item [{\rm (1)}] $K_0(s)$ is $LO$ for all $s \in \mathbb Q$ and $r$ is not the cabling slope of $K$ if $K$ is cabled, or
\item [{\rm (2)}] $P(U)(s)$ is $LO$ for all $s \in \mathbb Q$.
\end{enumerate}
\end{thm} 

Theorem \ref{thm: all implies all} is a consequence of Theorem \ref{thm: w = 0} and the following two propositions.

\begin{prop} 
\label{prop: $LO$ pattern}
Let $K = P(K_0)$ be a satellite knot and fix $r \in \mathbb Q$. Then $K(r)$ is $LO$ if $P(U)(r)$ is. 
\end{prop}

\begin{proof}
First, if $K$ is a cable knot, then by Lemma \ref{lemma: cabled implies never 1}(1) the pattern $P$ must be a cabled pattern. Since $\pi_1(P(U)(r))$ is left-orderable, then $r$ cannot be the cabling slope. So $K(r)$ is irreducible by Corollary \ref{cor: sch red sat}. Therefore, by \cite[Theorem 1.1]{BRW05}, it has a left-orderable fundamental group if and only if there is a non-trivial homomorphism of $\pi_1(K(r))$ to a left-orderable group. 

Let $F_0$ be the Seifert surface of $K_0$. We can collapse the exterior $X(K_0)$ of $K_0$ to a solid torus by collapsing a regular neighborhood of $\partial X(K_0)\cup F_0$ to a solid torus with an interior point removed and mapping the remaining part of $X(K_0)$ to the removed point. This shows that there is a slope-preserving degree one map $(X(K), \partial X(K)) \to (X(P(U)), \partial X(P(U)))$, and hence it induces a degree $1$ map $K(r) \to P(U)(r)$. The claim then follows. 
\end{proof}

\begin{prop} 
\label{prop: companion}
Let $K$ be a satellite knot with a pattern of winding number $w\geq 1$. Fix $r \in \mathbb Q$ and suppose that $r$ is not the cabling slope if $K$ is a cable knot. Then $K(r)$ is $LO$ if $K_0(r/w^2)$ is. 
\end{prop}

\begin{proof}
Suppose that $w \geq 1$ and $K_0(r/w^2)$ is $LO$. Then $r/w^2$ is $LO$-detected in $X_0$ by \cite[Corollary 8.3]{BC2}. Also, $\lambda_{P(r)} = r/w^2$ by (2.2.6), so if $T = \partial X_0$ is incompressible in $K(r)$, then $K(r)$ is $LO$ by Lemma \ref{lem: r/w^2 detected}. On the other hand, since $r$ is not the cabling slope, the case in Theorem \ref{thm: C(T)}(2)(b) doesn't occur.  Hence, if $T$ compresses in $K(r)$, then by Theorem \ref{thm: C(T)}(1) and (2)(c) we have $K(r) \cong K_0(r/w^2)$ and therefore is $LO$.
\end{proof}

We remark that if $T$ is incompressible in $K(r)$ then we only need $r/w^2$ to be $LO$-detected in $X_0$ to deduce that $K(r)$ is $LO$ by Lemma \ref{lem: r/w^2 detected}. 

\begin{proof}[Proof of Theorem \ref{thm: all implies all}]
If $P(U)(s)$ is $LO$ for all $s\in\mathbb{Q}$, then the conclusion follows from Proposition \ref{prop: $LO$ pattern}. Assume then that $K_0(s)$ is $LO$ for all $s\in \mathbb{Q}$. Note that if $K$ admits a pattern with winding number $w=0$, then $K(s)$ is $LO$ for any non-cabling slope $s\in \mathbb{Q}$ by Theorem \ref{thm: w = 0}. So we assume that $w\geq 1$. In this case, the conclusion follows from Proposition \ref{prop: companion}. 
\end{proof}

We remark that the method of proof of Proposition \ref{prop: $LO$ pattern} does not immediately apply to the $NLS$ and $CTF$ cases due to the following unanswered question.

\begin{question}
\label{que: degree nonzero map}
Let $M$ be a closed, connected, orientable, irreducible $3$-manifold. For $\ast \in \{NLS, \, CTF\}$, if there exists a nonzero degree map from $M$ to a closed, connected, orientable, irreducible $3$-manifold $N$ that has property $\ast$, does $M$ have property $*$? 
\end{question}

Since Question \ref{que: degree nonzero map} is known to have a positive answer when $* = LO$ \cite[Theorem 3.7]{BRW05}, the $L$-space conjecture predicts positive answers for both $NLS$ and $CTF$, though currently there are very few results in this direction. See \cite{LM18,HLL22} for partial answers to Question \ref{que: degree nonzero map} in the $NLS$ case.

The argument of Proposition \ref{prop: companion} can be applied to prove the analogous statement for the $NLS$ property. For $CTF$, we need to replace the condition that $K(r/w^2)$ is $CTF$ by  the a priori stronger condition that $r/w^2$ is a strongly $CTF$-detected slope. See Proposition \ref{prop: ctf companion}.

\subsection{\texorpdfstring{$LO$}{LO} and \texorpdfstring{$NLS$}{NLS} surgery on satellite knots and JSJ-graphs}
\label{subsec: JSJ graph}

\begin{thm} 
\label{thm: split} 
Suppose that the exterior $X(K)$ of a knot $K$ contains disjoint essential tori $T_1, T_2$ which together with $\partial X(K)$ cobound a connected submanifold of $X(K)$. Let $r \in \mathbb Q$ and assume that $r$ is not the cabling slope if $K$ is a cable knot. Then $K(r)$ is $LO$ and $NLS$. 
\end{thm}
 
Before proving this theorem, we deduce an application. 

\begin{figure}[ht]
\centering
\includegraphics[scale=1]{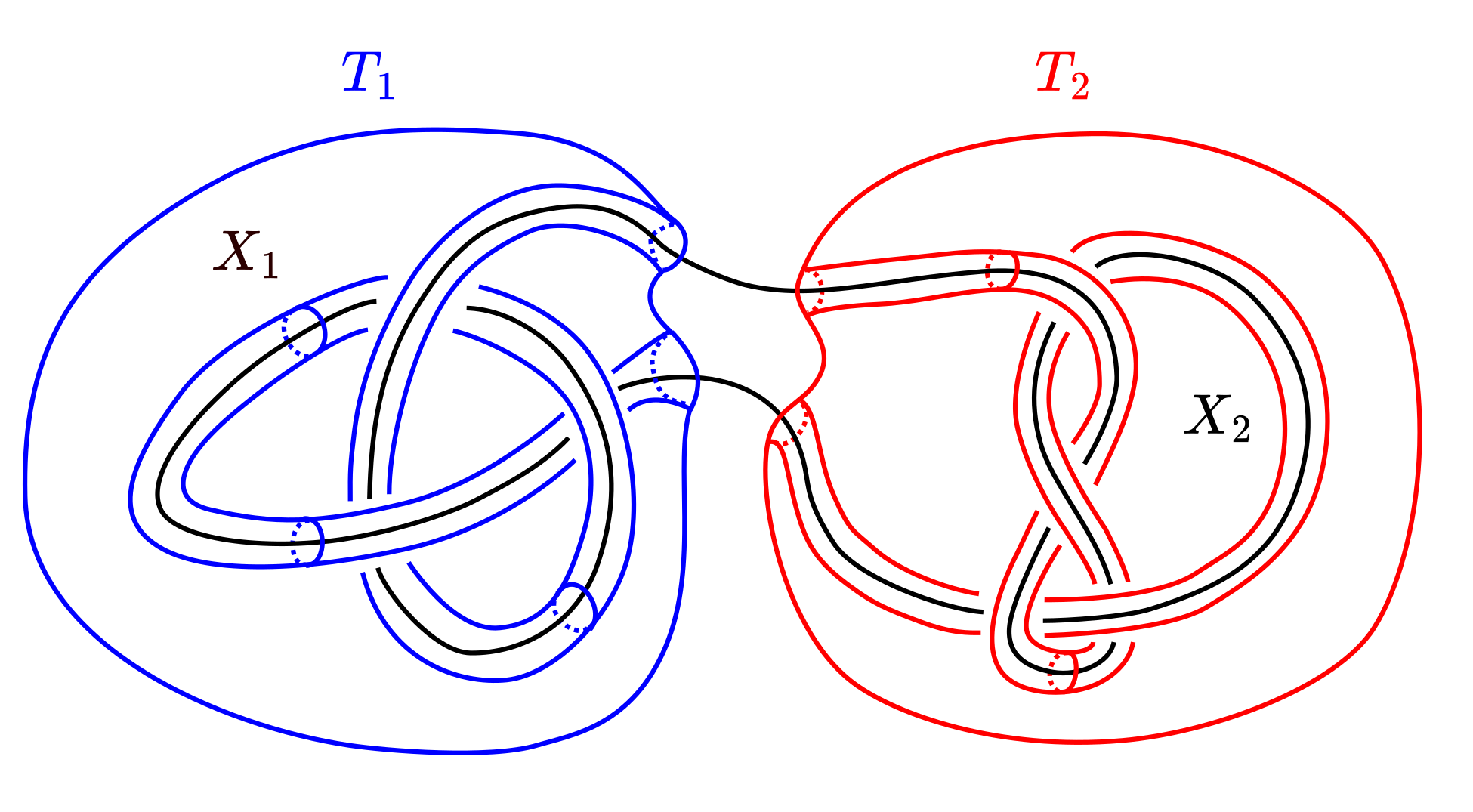}
\caption{$T_1$ and $T_2$ are two ``swallow-follow'' tori in the exterior of $K = K_1 \# K_2$. Though $\partial X(K)$ is not drawn in the figure, it is obvious that $T_1$, $T_2$ and $\partial X(K)$ cobound a submanifold in $X(K)$. In this case, the submanifold is homeomorphic to the product of a pair of pants and $S^1$.}
\label{fig: composite knots}
\end{figure}

Figure \ref{fig: composite knots} illustrates that the exterior $X(K)$ of the connected sum $K$ of two non-trivial knots contains two essential swallow-follow tori $T_1, T_2$ which together with $\partial X(K)$ cobound a connected submanifold of $X(K)$. Note, moreover, that these tori have winding number $1$, so $K$ is not a cable knot. Hence as a special case of Theorem \ref{thm: split} we deduce that all rational surgeries on $K$ are $NLS$ and $LO$. This gives another proof of Krcatovich's result that $L$-space knots are prime \cite{Krc15} and its $LO$ counterpart. In summary, we have the following corollary from Theorem \ref{thm: split}.

\begin{cor}
\label{cor: nlo irred prime}
All rational surgeries on a composite knot are $LO$ and $NLS$. 
\end{cor}

\begin{proof}[Proof of Theorem \ref{thm: split}]
We can assume that $X(K)$ contains no winding number $0$ essential tori by Theorem \ref{thm: w = 0}. Since $K(0)$ is irreducible \cite{Gab3}, $K(0)$ is $LO$ by \cite[Theorem 1.1]{BRW05}. And $K(0)$ is not an $L$-space by definition. Hence, we can assume that $r \ne 0$. 

Write $X(K) = Y \cup_{T_1} X_1 \cup_{T_2} X_2$, where $Y$ is the connected submanifold cobounded by $T_1$, $T_2$, and $\partial X(K)$ and $X_1, X_2$ are non-trivial knot exteriors with $\partial X_i = T_i$. By Corollary \ref{cor: sch red sat}, $K(r) = X_1 \cup_{T_1} Y(r) \cup_{T_2} X_2$ is irreducible,  where $Y(r)$ is the $r$-Dehn filling of $Y$ along $\partial X(K)$. Note that $T_i$ bounds a solid torus in $S^3$, which we denote by $V_i$, and $K\subset V_i$, $i=1,2$.

We claim that $Y(r)$ is irreducible. Otherwise, the irreducibility of $K(r)$ implies that there exists a $2$-sphere in $Y(r)$ that bounds a $3$-ball $B$ in $K(r)$ but not in $Y(r)$. Then $B$ must contain some $T_i$, say $T_1$, and therefore $T_1$ is compressible in $K(r)$. Let $P_1$ denote $K$ viewed as a knot in $V_1$. Then $K(r) = X_1\cup_{T_1} P_1(r)$, where $P_1(r) = Y(r)\cup_{T_2} X_2$. Since $T_1$ is incompressible in $X_1$, it is compressible in $P_1(r)$ so either case (a) or case (b) of Theorem \ref {thm: sch red} applies to $P_1$ in $V_1$. Case (b) doesn't hold by our hypothesis that $r$ is not a cabling slope. Hence Theorem \ref{thm: sch red}(a) holds and so $P_1$ is braided in $V_1$. But then by \cite[Lemma 3.1]{Be91}, the torus $T_2$ bounds a solid torus in $V_1$, contradicting the fact that it bounds $X_2$.

By a result of Schubert, see \cite[Proposition 2.1]{Bu06}, we can find meridional disks $D_i$ of the solid tori $V_i = S^3\setminus \text{int}(X_i) $  which are contained in $S^3 \setminus (\text{int}(X_1) \cup \text{int}(X_2)) = Y(\mu)$, where $\mu$ is the meridional class of $K$. This implies that the meridional class $\mu_1$ of $X_1$, which is represented by $\partial D_1$, is homologous to $w_1 \mu$ in $Y$, where $w_1 \ne 0$ is the winding number of $T_1$ in $X(K)$. Similarly the meridional class $\mu_2$ of $X_2$ is homologous to $w_2 \mu$ in $Y$, where $w_2 \ne 0$. Thus $w_2 [\mu_1] = w_1 [\mu_2] {\text{ in }} H_1(Y)$,and hence also in $H_1(Y(r))$.

Let $U_i$ be the twisted $I$-bundle over the Klein bottle, $\lambda_{U_i}$ the longitudinal slope of $U_i$, and let $Z$ be the manifold $U_1\cup_{T_1}Y(r)\cup_{T_2} U_2$ obtained by gluing $U_i$ to $Y(r)$ along $T_i$ with $\lambda_{U_i}$ identified with $\mu_i$. Then $Z$ is irreducible. Also, $H_1(Z; \bQ) \cong H_1(Y(r); \bQ)/ \langle[\mu_1],[\mu_2]\rangle$. Since $r \ne 0$, $Y(r)$ is a $\bQ$-homology cobordism between $T_1$ and $T_2$, so $\dim H_1(Y(r); \bQ) = 2$. Therefore, since as noted above $[\mu_1]$ and $[\mu_2]$ are linearly dependent in $H_1(Y(r); \bQ)$, $\dim H_1(Z; \bQ) \ge 1$. Hence $Z$ is $LO$ and $NLS$. Then by \cite[Definition 4.6, Proposition 6.13]{BGH21}, the multislope $([\mu_1], [\mu_2])$ is $LO$- and $NLS$-detected in $Y(r)$. On the other hand, each $[\mu_i]$ is $LO$- and $NLS$-detected in $X_i$ by Theorem \ref{thm: meridional detn}. Therefore, using the multislope gluing theorem \cite[Theorem 7.6]{BGH21}, we have $K(r) = X_1 \cup_{T_1} Y(r) \cup_{T_2} X_2$ is $LO$ and $NLS$. This completes the proof. 
\end{proof}

Recall the rooted JSJ graph of a knot described in the introduction. If the exterior $X(K)$ of a satellite knot $K$ contains disjoint essential tori $T_1, T_2$ which together with $\partial X(K)$ cobound a connected submanifold $Y$ of $X(K)$, then the rooted JSJ graph of $K$ cannot be a rooted interval. See Figure \ref{fig: jsj graph} below. Conversely, if the JSJ graph of $K$ is not a rooted interval, it is easy to see that there must exist $T_1$ and $T_2$ that together with $\partial X(K)$ cobound a connected submanifold $Y$. 
\begin{figure}[ht]
\centering 
\includegraphics[scale=1.3]{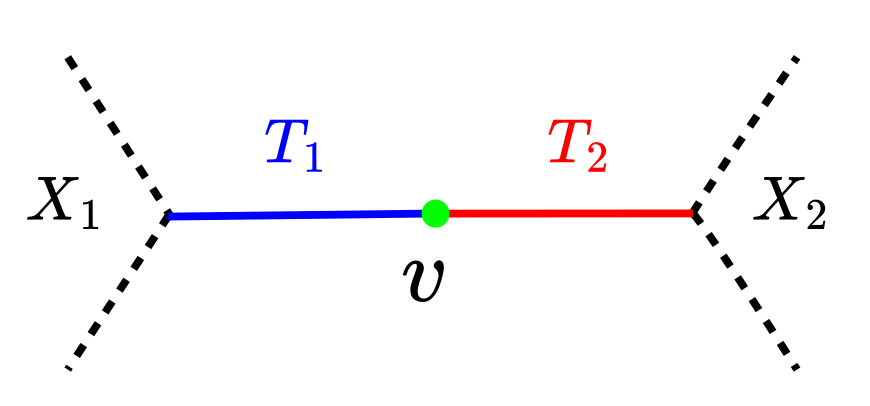}
\caption{The figure illustrates the case that $Y$ is atoroidal and the green vertex $v$, representing the submanifold $Y$, is the root of the tree. In general, there could be other branches at the vertex $v$ in addition to the two branches depicted in the figure. The root of the tree is contained in one of these additional branches. The submanifold $Y$ is just $X(K)\setminus {\text{int}(X_1)\cup \text{int}(X_2)}$, corresponding to the union of all these additional branches.}
\label{fig: jsj graph}
\end{figure}

Here are two immediate consequences of Theorem \ref{thm: split}. 

\begin{thm}
\label{thm: jsj graph interval $L$-space knot}
Suppose that $K$ is a satellite $L$-space knot. Then the JSJ graph of $K$ is  a rooted interval.
\end{thm}

\begin{thm}
\label{thm: jsj lo}
Suppose that $K$ is a satellite knot which admits an irreducible rational surgery that is not $LO$. Then the JSJ graph of $K$ is a rooted interval.
\end{thm}

\section{Co-oriented taut foliations and surgeries on satellite knots} 
\label{sec: surgery on satellite knots ctf}

The analogous behaviour of detection and gluing in the $LO, NLS$, and $CTF$ contexts (cf. Theorems \ref{thm: * gluing} and \ref{thm: meridional detn}) means that the arguments used to verify the properties $LO$ and $NLS$ in the proofs of Theorems \ref{thm: w = 0} and \ref{thm: split} can be used to deduce comparable results for the property $CTF$ as long as  the meridional slope of the companion knots which arise are $CTF$-detected. In fact, we conjecture that the meridional slope of any non-trivial knot in the $3$-sphere is $CTF$-detected (cf. Conjecture \ref{conj: CTF meridional detection}) and verified this for fibred knots (see Theorem \ref{thm: meridional detn}). Consequently we deduce the following three results.  

\begin{thm}
\label{thm: ctf winding zero fibered}
Suppose that $K$ is a satellite knot whose exterior contains a winding number zero essential torus which bounds the exterior of a fibred knot in $X(K)$. Then $K(r)$ is $CTF$ unless $K$ is a cable knot and $r$ is the cabling slope.
\end{thm} 

\begin{proof}
The proof is identical to that of Theorem \ref{thm: w = 0}. The fibring assumption is required so that we can apply Theorem \ref{thm: meridional detn}.
\end{proof}

\begin{thm} 
\label{thm: ctf split fibered} 
Suppose that the exterior $X(K)$ of a knot $K$ contains disjoint essential tori $T_1, T_2$ which together with $\partial X(K)$ cobound a connected submanifold of $X(K)$, and each of which bounds the exterior of a fibred knot. Then $K(r)$ is $CTF$ unless $K$ is a cable knot and $r$ is the cabling slope.
\end{thm}

\begin{proof}
The proof is identical to that of Theorem \ref{thm: split}.  
\end{proof}

An immediate corollary of Theorem \ref{thm: ctf split fibered} is the following $CTF$ analogue of Theorem \ref{thm: jsj lo}.

\begin{thm}
\label{thm: ctf rooted interval}
Suppose that $K$ is a fibered satellite knot which admits an irreducible rational surgery which is not $CTF$. Then the JSJ graph of $K$ is a rooted interval.
\end{thm} 

It is interesting to note that since $CTF$ implies $NLS$, the only additional facts needed to deduce Theorem \ref{thm: jsj graph interval $L$-space knot} from Theorem \ref{thm: ctf rooted interval} are that an $L$-space knot is fibered and has at least one irreducible $L$-space rational surgery. 

Theorem \ref{thm: ctf rooted interval} has the following corollary, which was previously proven by Delman and Roberts in \cite{DR21} using different methods.

\begin{cor} {\rm (Delman-Roberts \cite[Corollary 6.2]{DR21})}
\label{cor: ctf composite fibered}
All rational Dehn surgeries on a composite fibred knot are $CTF$. 
\end{cor}

Next, we consider the $CTF$ analogue of Proposition \ref{prop: companion}. The argument used in the proof of that proposition cannot be applied to the $CTF$ case owing to the fact that unlike the $NLS$ and $LO$ cases, it is unknown whether $K(r)$ being $CTF$ implies that $r$ is $CTF$-detected (see Remark \ref{rem: detection and Dehn fillings}). We can apply the proof if we use a strengthened form of $CTF$-detection though.  

In \cite{BGH21}, we defined a slope $\alpha$ on the boundary of a knot manifold $M$ to be {\it strongly $CTF$-detected} if there exists a co-orientable taut foliation on $M$ that  transversely intersects $\partial M$ in a linear foliation by simple closed curves of slope $\alpha$. If $\alpha$ is strongly $CTF$-detected, it is clear that it is $CTF$-detected and, moreover, that $M(\alpha)$ is $CTF$. It is unknown whether $M(\alpha)$ being $CTF$ implies that $\alpha$ is strongly $CTF$-detected (cf. Remark \ref{rem: detection and Dehn fillings}).  

\begin{prop}
\label{prop: ctf companion}
Suppose that $K = P(K_0)$ is a satellite knot where $P$ has winding number $w\geq 1$. Then $K(r)$ is $CTF$ if $r/w^2$ is a strongly $CTF$-detected slope of $X_0$ and $r$ is not the cabling slope when $K$ is a cable knot.
\end{prop}
 
\begin{proof}
Let $K(r) = P(r) \cup_{T_0} X_0$ as before. Note that the longitudinal slope $\lambda_{P(r)}$ of $P(r)$ is the slope $r/w^2$ on $\partial X_0$ by (\ref{eqn: slope as a rational number}). If $T_0 = \partial X_0$ is incompressible in $K(r)$, then $K(r)$ is $CTF$ by Lemma \ref{lem: r/w^2 detected}. Suppose, on the other hand, that $T_0$ is compressible in $K(r)$. Since we have assumed that $r$ is not a cabling slope, Theorem \ref{thm: C(T)}(1) and (2)(c) imply that $K(r) \cong K_0(r/w^2)$ and therefore is $CTF$ by the hypothesis that $r/w^2$ is strongly $CTF$-detected. 
\end{proof}

\begin{remark}
\label{remark: incomp torus}
The proof of Proposition \ref{prop: ctf companion} shows that if $\partial X_0$ is incompressible in $K(r)$, we need only assume that $r/w^2$ is $CTF$-detected in $X_0$ to conclude that $K(r)$ is $CTF$. 
\end{remark}

Delman and Roberts call a knot $K$ {\em persistently foliar} if all rational slopes of $K$ are strongly $CTF$-detected. The theorem below is an immediate consequence of Proposition \ref{prop: ctf companion}.

\begin{thm} 
\label{thm: satellite ctf}
Suppose that $K$ is a satellite knot with a persistently foliar companion whose associated pattern has winding number $w\geq 1$. Then $K(r)$ is $CTF$ for any $r\in \mathbb{Q}$ which is not the cabling slope if $K$ is a cable knot. 
\end{thm}

In \cite{DR21}, Delman and Roberts have shown that if $K$ is persistently foliar, then the same is true for the connected sum of $K$ with any other knot. Consequently, every rational surgery on such a knot admits a co-oriented taut foliation. One can deduce the same conclusion using Theorem \ref{thm: satellite ctf}. 

\begin{cor}{\rm (cf. Delman-Roberts \cite[Corollary 4.2]{DR21})}
\label{cor: dr composite}
Each rational surgery on a composite knot with a persistently foliar summand is $CTF$. 
\end{cor}

\begin{proof}
Let $K$ be a composite knot with a persistently foliar summand $K_0$. View $K$ as a satellite knot with companion $K_0$. In this case the associated winding number is $1$, so $K$ is not a cable knot (cf. Lemma \ref{lemma: cabled implies never 1}(3)). The desired conclusion now follows from Theorem \ref{thm: satellite ctf}.
\end{proof}

Lastly, we use Proposition \ref{prop: ctf companion} to deduce Proposition \ref{prop: ctf fibred} below, which is based on a result of Roberts but also strengthens it. More precisely, in \cite[Theorem 4.7]{RobertsSurfacebundle2}, Roberts shows that if $K$ is a fibred  knot whose monodromy has positive (resp. negative, resp. zero) fractional Dehn twist coefficient, then any slope $r\in (-\infty, 1)$ (resp. $(-1, \infty)$, resp. $(-\infty, \infty)$) is strongly $CTF$-detected in $X(K)$\footnote{Though the statement of Theorem 4.7 in \cite{RobertsSurfacebundle2} requires the knot to be hyperbolic, the argument in \cite{RobertsSurfacebundle2} holds for all non-trivial fibred knots (cf. \cite[Theorem 1.4]{DR21}).}. This, together with Proposition \ref{prop: ctf companion},  implies the following. 

\begin{prop}
\label{prop: ctf fibred}
Let $K = P(K_0)$ where $K_0$ is fibered and $P$ has winding number $w$. Let $c$ be the fractional Dehn twist coefficient of the monodromy of $K_0$. Suppose that $r$ is not the cabling slope if $K$ is a cable knot. Then $K(r)$ is $CTF$ if either
\begin{enumerate}[leftmargin=*] 
\setlength\itemsep{0.3em} 
\item[{\rm (1)}] $c > 0$ and $r \in (-\infty, w^2]$, or
\item[{\rm (2)}]$c < 0$ and $r \in [-w^2, \infty)$, or
\item[{\rm (3)}] $c = 0$ and $r \in (-\infty, \infty)$.
\end{enumerate}

\end{prop}

\begin{proof}
If $w=0$, then Theorem \ref{thm: ctf winding zero fibered} shows that $K(r)$ is $CTF$ for any rational slope $r$ that is not a cabling slope. So we assume that $w\geq 1$.

If $ c > 0$ then by \cite[Theorem 4.7]{RobertsSurfacebundle2}, any slope in $(-\infty, 1)$ is strongly $CTF$-detected in $X_0$. Then $K(r)$ is $CTF$ for any slope $r\in (-\infty, w^2)$ that is not a cabling slope by Proposition \ref{prop: ctf companion}. 

To deal with the case $r = w^2$ we use the facts that $\partial X_0$ is incompressible in $K(w^2)$ by Lemma \ref{lemma: incomp torus}, and that the slope 1 is $CTF$-detected in $X_0$ by Theorem \ref{thm: meridional detn}. Then $K(w^2)$ is $CTF$ by Remark \ref{remark: incomp torus}. 

This proves (1). (2) is completely analogous, and (3) follows in the same way using the fact that if $c = 0$ then all rational slopes on $\partial X_0$ are strongly $CTF$-detected in $X_0$ \cite[Theorem 4.7]{RobertsSurfacebundle2}. 
\end{proof}

We note that fibred strongly quasipositive satellite knots satisfy the hypothesis (1). For, suppose $K=P(K_0)$ is fibred and strongly quasipositive. Then $K_0$ is fibred, and it is strongly quasipositive by the main result of \cite{Rudolph92}. Its monodromy is right-veering by \cite{Hed10,HKMI} and hence the fractional Dehn twist coefficient of its monodromy is positive.

It is known that positive $L$-space knots (i.e. knots with a positive $L$-space surgery) are fibred and strongly quasipositive \cite{Ni07, Hed10}. We say that a satellite knot is a {\it positive satellite $L$-space knot} if it is also a positive $L$-space knot. Then Proposition \ref{prop: ctf fibred}(1) applies to positive satellite $L$-space knots. In fact, Lemmas \ref{lemma: cab $L$-space knot} and \ref{lemma: $L$-space knot $w^2$} below show that in that case, the possibility that $r$ is a cabling slope does not arise. Hence we have the following corollary, which by Lemma \ref{lemma: $L$-space knot $w^2$} is predicted by the $L$-space Conjecture.

\begin{cor}
 \label{cor: $L$-space ctf}   
If $K$ is a positive satellite $L$-space knot with pattern of winding number $w$ then $K(r)$ is $CTF$ for each rational $r \in (-\infty, w^2]$.
\end{cor}

\begin{lemma}
\label{lemma: cab $L$-space knot}
  Let $K$ be a positive satellite $L$-space cable knot. Then surgery on $K$ along the cabling slope is an $L$-space. 
\end{lemma}

\begin{proof}
We have $K = C_{m,n}(K_0)$. Let $g_0 = g(K_0)$. By \cite{Hom11},  the assumption that $K$ is a positive $L$-space knot implies that $n/m > 2g_0 - 1 \ge 1$, so $m + n > 2mg_0$.

By \cite{Schu53}, $$g(K) = g(C_{m,n}(U)) + mg_0$$

Since $C_{m,n}(U)$ is the torus knot $T(m,n)$ we get 
\begin{align*}
2g(K) - 1 & = (m-1)(n-1) + 2mg_0 -1 \\
  &= mn - (m+n-2mg_0) < mn
\end{align*}
\end{proof}

\begin{lemma}
\label{lemma: $L$-space knot $w^2$}
Let $K$ be a positive satellite $L$-space knot with pattern of winding number $w$. Then $w^2 < 2g(K) - 1$.
\end{lemma}

\begin{proof}
We have $K = P(K_0)$, where $P$ has winding number $w$. Let $g_0 = g(K_0)$. 

By \cite{Schu53}, $$g(K) = g(P(U)) + wg_0$$  

Hence $$2g(K) - 1 = 2g(P(U)) - 1 +2wg_0$$

By \cite[Theorem 7.3(i)]{BMot19}, 
\begin{align*}
    2g(P(U)) - 1 & > w(w - 1)(2g_0 - 1) - w \\
    & = (w-1)2wg_0 - w^2
\end{align*}                                     
Hence $2g(K) - 1 > w^2(2g_0 - 1) \ge w^2$.
\end{proof} 

Lastly, we show that surgery on a positive satellite $L$-space knot with surgery coefficient at most 9 is $CTF$ if and only if it is $NLS$ (Corollary \ref{cor: ctf nls r<=9}). This is an immediate consequence of Theorem \ref{thm: sat $L$-space ctf} below.  

\begin{thm}
\label{thm: sat $L$-space ctf}
   Let $K= P(K_0)$ be a positive satellite $L$-space knot. Then $K(r)$ is $CTF$ if $r \in (-\infty,9]$ unless $K(r)$ is an $L$-space. The latter happens exactly when $K = C_{2,n}(T(2,3))$ with $n = 3,5,$ or $7$ and $r \in [n+2,9]$.
\end{thm}

\begin{proof}
Let $w$ denote the winding number of $P$. When $w \ge 3$, the result follows from Corollary \ref{cor: $L$-space ctf}. So we can assume that $w = 2$ by Corollary \ref{cor: w>1 nls}. 

Since $K$ is a positive $L$-space knot, the companion $K_0$ is a positive $L$-space knot by \cite[Theorem 1.15]{HRW1}, and $P$ is braided in its solid torus by \cite[Theorem 1.17]{BMot19} (also see footnote $6$). So $P = C_{2,n}$ for some odd $n$.  We consider three cases: $g(K_0) = 1$, $g(K_0) = 2$ and $g(K_0) >2$. 

  \begin{case} 
  \label{case: g = 1}
 $g(K_0) = 1$. 
  \end{case}
Then $K_0 = T(2,3)$ by \cite{Gh08}. Since $K$ is a positive $L$-space knot, \cite{Hom11} implies that $n/2 > 2g(K_0) -1 = 1$, so $n \ge 3$. Also, $C_{2,n}(U) = T(2,n)$, and therefore $g(K) = g(C_{2,n}(U)) + 2g(K_0) = (n-1)/2 + 2 = (n+3)/2$ by \cite{Schu53}. If $r \neq 2n$, then $K(r)$ is either Seifert fibred  or a graph manifold, therefore it is $CTF$ if and only if it is $NLS$ (\cite{LS07, BGW13, BC17, HRRW15}), and the latter happens if and only if $r < 2g(K) - 1 = n + 2$ by \cite{OS05-lens,OS11}. For $n + 2 \le r \le 9$, we obtain the exceptions listed in the theorem. When $r=2n \leq 9$, we have $n=3$ and $r=2n=6$. So $K(6) = L(2,3)\#T(2,3)(3/2)$ is an $L$-space and is not $CTF$.
\begin{case}
$g(K_0) = 2$. 
 \end{case}
 Then $K_0 = T(2,5)$ by \cite{FRW22}. By \cite{Hom11}, $n/2 > 2g(K_0) -1 = 3$, so $n \ge 7$ and the cabling slope $2n\geq 14$. Also, $g(K) = (n-1)/2 + 2 \cdot 2 = (n+7)/2$, and hence, as in the case above, $K(r)$ is $CTF$ if and only if $r < 2g(K) - 1 = n+6$ by \cite{OS05-lens,OS11, LS07, BGW13, BC17, HRRW15}. Therefore $K(r)$ is $CTF$ for $r < 13$.
  
 \begin{case}
    $g(K_0) >2$. 
     \end{case}
In this case, we will show that $K(r)$ is $CTF$ for any $r < 13$.
By \cite{Hom11}, $n/2 > 2 \cdot 3 -1 = 5$, so $n \ge 11$ and the cabling slope is $2n\geq 22$. 
We have $K(r) = C_{2,n}(r) \cup_{T} X_0$, where $X_0$ is the exterior of $K_0$.

We first claim that $T$ must be incompressible in $K(r)$ when $r < 21$. To see this, set $r=p/q$ and note that if $T$ compresses in $K(r)$, we have $|p - 2nq| = \Delta(p/q,2n/1) \leq 1$ and therefore $|p/q - 2n| \leq 1/q \leq 1$ by Theorem \ref{thm: C(T)}(2)(a). Since $2n \geq 22$, $r = p/q \geq 21$, which proves the claim. 

 We will show that when $r < 13$ there exists a slope $s$ on $T$ that is $CTF$-detected in both $C_{2,n}(r)$ and $X_0$. Since $T$ is incompressible in $K(r)$, the claim follows from Theorem \ref{thm: * gluing}. 

Since $K_0$ is a positive $L$-space knot, it is fibred with right-veering monodromy by \cite{Ni07, HKMI, Hed10}. So by \cite{RobertsSurfacebundle2}, all slopes in $(-\infty, 1)$ are strongly $CTF$-detected in $X_0$, and hence $CTF$-detected in $X_0$. In addition, slopes $\mu$ and $1$ are $CTF$-detected by Theorem \ref{thm: meridional detn}. Therefore, the set of $CTF$-detected slopes of $X_0$ contains $[- \infty,1]$. 

To analyse the set of $CTF$-detected slopes in $C_{2,n}(r)$, consider $K' = C_{2,n}(T(2,3))$. Then $K'(r) = C_{2,n}(r) \cup_{T'} X(T(2,3))$. Since $n/2 \ge 11/2 > 2g(T(2,3)) - 1  = 1$, \cite{Hed09} implies that $K'$ is an $L$-space knot, and as we argued in Case \ref{case: g = 1},  $K'(r)$ is $NLS$ if and only if $r < n + 2$. Since $n \geq 11$ in this case, and $K'(r)$ is $NLS$ for $r < 13$. Therefore, by \cite[Theorem 1.14]{HRW1}, for $r < 13$, there is a slope $s$ on $T'$ that is $NLS$-detected in both $C_{2,n}(r)$ and $X(T(2,3))$. Since $C_{2,n}(r)$ and $X(T(2,3))$ are Seifert fibre spaces, $s$ is $CTF$-detected in both $C_{2,n}(r)$ and $X(T(2,3))$ \cite[Theorem 1.6]{BC17}. Now the $CTF$-detected slopes for $X(T(2,3))$ are precisely those in the interval $[- \infty, 1]$ (with respect to the standard meridional and longitudinal basis). So the slope $s$ as a slope on $\partial X_0 = T = \partial C_{2,n}(r) = T'$ is also contained in $[-\infty, 1]$. 

Therefore there is a slope $s$ on $T$ that is $CTF$-detected in both $C_{2,n}(r)$ and $X_0$. Hence $K(r)$ is $CTF$ for $r < 13$. 
\end{proof}

\begin{cor}
\label{cor: ctf nls r<=9}
If $K$ is a positive satellite $L$-space knot then for any $r \in (-\infty, 9]$, $K(r)$ is $CTF$ if and only if it is $NLS$. 
\end{cor}

\section{Left-orderable \texorpdfstring{$p/q$}{p/q}-surgeries on knots with \texorpdfstring{$p$}{p} small}
\label{sec: $LO$ surgery for p small} 

Left-orderable surgery on torus knots is completely understood. In fact, given a coprime pair $m, n \geq 2$, the $(m, n)$ torus knot $T(m,n)$ is an $L$-space knot. As the $L$-space conjecture is known for Seifert manifolds (\cite{BRW05,LS07,BGW13}), if $p/q \in \mathbb Q$ then $T(m, n)(p/q)$ is $LO$ if and only if $p/q < 2g(T(m, n)) -1 = mn - (m+n)$.

In this section, we consider $p/q$-surgeries on hyperbolic and satellite knots with $p$ small and $q \ne 0$. Theorem \ref{thm: small p} is a simplified version of our main result, Theorem \ref{thm: $LO$ surgery small p full}.

\begin{thm}
\label{thm: $LO$ surgery small p full}
Let $K$ be a non-trivial knot and let $E$ denote the set of rational numbers $p/q$ such that $p = 1$ or $2$, $q \ne 0$, and $K(p/q)$ is not $LO$. If $E$ is non-empty, then one of the following occurs: 
\vspace{-.2cm} 
\begin{enumerate}[leftmargin=*]
\setlength\itemsep{0.3em}
    \item[{\rm (1)}] $K = T(2,3\varepsilon)$ for some $\varepsilon  = \pm1$ and  $E = \{ \varepsilon, 2\varepsilon\};$ or
    \item[{\rm (2)}] $K$ is a $(2,\varepsilon)$-cable for some $\varepsilon = \pm1$ and $E = \{2\varepsilon\};$ or 
    \item[{\rm (3)}]  $K$ is hyperbolic and there is some $\varepsilon \in \{\pm 1\}$ such that $E$ is contained in either $\{\varepsilon, 2\varepsilon\}$ or $\{\varepsilon/2, 2\varepsilon/3, \varepsilon, 2\varepsilon\}$.
\end{enumerate} 
\end{thm}

\begin{remark}
\label{rem: mostly nec}
Since a hyperbolic $L$-space knot has genus at least $3$ by \cite{Gh08, FRW22}, the $L$-space Conjecture predicts that for any hyperbolic knot $K$, $K(r)$ is $LO$ either for all $r \in (-\infty, 5)$ or for all $r \in (-5, \infty)$. Thus it is expected that the potential exceptional slopes listed in case (3) of the theorem do not arise.  
\end{remark}

\subsection{\texorpdfstring{$LO$}{LO} \texorpdfstring{$p/q$}{p/q}-surgeries on satellite knots for \texorpdfstring{$p$}{p} small}
\label{subsec: $LO$ surgery on satellites for p small}

In \cite{BGH21}, we proved the following theorem. 

\begin{thm}[Theorem 2.1 in \cite{BGH21}]
    \label{thm: $LO$ small HRRW15}
Suppose that $W$ is a closed, connected, orientable, irreducible, toroidal $3$-manifold. If $|H_1(W)|\leq 4$, then $W$ is LO. 
\end{thm} 

It follows that if $1 \leq p \leq 4$ and $K(p/q)$ is irreducible and toroidal, then $K(p/q)$ is $LO$. The following result is then a consequence of the discussion of irreducible and toroidal surgery on satellite knots in \S \ref{sec: background}.

\begin{prop}
\label{prop: satellite p small}
Let $K$ be a satellite knot, and suppose that $p/q$ is a reduced fraction with $1 \leq p \leq 4$ and $|q| \geq 2$.
Then $K(p/q)$ is $LO$ unless there is an $\varepsilon \in \{\pm 1\}$ such that $K$ is the $(2, \varepsilon)$-cable of a hyperbolic knot $K_0$, where $p/q = 3\varepsilon/2$ and $K_0(3\varepsilon/8)$ is not $LO$.
\end{prop}

\begin{proof}
Since $1 \leq p \leq 4$, Theorem \ref{thm: $LO$ small HRRW15} implies that $K(p/q)$ will be $LO$ as long as $K(p/q)$ is irreducible and toroidal. The former is guaranteed by the condition that $|q| \geq 2$ by Corollary \ref{cor: sch red sat}, so we need only analyse when the latter occurs. Theorem \ref{thm: C(T)} shows that $K(p/q)$ is toroidal if $K$ is not a cable knot, so suppose that it is, say $K = C_{m,n}(K_0)$ where $K_0$ is non-trivial and $m \geq 2, n \neq 0$. 

Assume that $|q| \geq 2$ and that $K(p/q)$ is not $LO$, and hence is atoroidal.
Then $T$, the boundary of the exterior of $K_0$, must compress in $C_{m,n}(p/q)$. By Theorem \ref{thm: sch red}(2)(a), and since $|q| \ge 2$, this happens if and only if $qmn = p \pm 1$. Then 
$$4|n| \leq |q|m|n| = |p \pm1| \leq 5$$
and therefore $n = \varepsilon$, $q = 2\varepsilon$, $m = 2$ and $p = 3$, for some $\varepsilon \in \{\pm 1\}$. Hence
$$K(p/q) = K(3\varepsilon/2) \cong K_0(3\varepsilon/8)$$
Since $|3\varepsilon/8| < 1$, $3\varepsilon/8$-surgery on a torus knot is $LO$ and therefore $K_0$ must be hyperbolic.
\end{proof}

\begin{remarks}

$\;$

\begin{enumerate}[leftmargin=*]
\setlength\itemsep{0.3em}
\item Since $-1 < 3 \varepsilon /8 < 1$, $3\varepsilon/8$-surgery on any non-trivial knot is not an $L$-space. Hence in Proposition \ref{prop: satellite p small}, $K(p/q) = K(3\varepsilon/2) \cong K_0(3\varepsilon/8)$ is not an $L$-space and therefore the $L$-space conjecture predicts that it is $LO$. Thus the potential exception listed in the proposition is not expected to exist. 
\item Similar arguments can be used to show that if $K(p)$ is not $LO$ for some $1 \leq p \leq 4$, then 
\medskip
\begin{enumerate}
\item if $K$ is not a cable knot, there are at most two values of such $p$ for which $K(p)$ is not $LO$, and if two, they are successive integers; 
\item if $K = C_{m,n}(K_0)$, then $K(p)$ is $LO$ unless $|n| = 1$ and either 
\begin{enumerate}
    \setlength\itemsep{0.3em}
        \item[{\rm (i)}] $p$ is the cabling slope $mn$;
        \item[{\rm (ii)}]  $p = mn \pm 1$ and $K(p) = K(mn \pm 1) \cong K_0((mn \pm 1)/m^2)$, where $K_0$ is a hyperbolic knot, and $K_0((mn \pm 1)/m^2)$ is not $LO$. Since $\frac{|mn \pm 1|}{m^2} \leq \frac{m +1}{m^2} \leq 3/4$, we do not expect such a knot $K_0$ to exist.
    \end{enumerate}
\end{enumerate}

\medskip
\noindent

\end{enumerate}
\end{remarks}

\subsection{\texorpdfstring{$LO$}{LO} \texorpdfstring{$p/q$}{p/q}-surgery on hyperbolic knots when \texorpdfstring{$p\leq 2$}{p<=2}} 
\label{subsec: $LO$ surgery hyperbolic small p}
The question of which surgeries on a hyperbolic knot $K$ are $LO$ is wide open in general. In the case that $p = 1$ or $2$, $q \ne 0$, and $K$ is hyperbolic, it is known that $K(p/q)$ is never an $L$-space. Thus we expect all such surgeries to yield manifolds that are $LO$. In this section, we will present some results regarding the left-orderability of the fundamental groups of these manifolds and describe their proofs using known results on pseudo-Anosov flows. See \cite{Mosher96} and \cite[\S 6.6]{Cal07} for background information on pseudo-Anosov flows. 

A pseudo-Anosov flow $\Phi_0$ on the complement of a knot $K$ determines an essential lamination in the interior of $M$, i.e., the (un)-stable lamination of $\Phi_0$. The connected complementary component $C$ of the lamination that contains $\partial M$ is a bundle over $S^1$ whose fiber is an ideal polygon with an open disk removed from its interior. It is easy to see that this component is homeomorphic to $\partial M\times [0, 1]$ with a collection of simple closed curves removed on $\partial M\times \{1\}$, where $\partial M$ is identified with $\partial M\times \{0\}$. This collection of simple closed curves is called the {\it degeneracy locus} of the pseudo-Anosov flow. See \cite[Definition 6.42]{Cal07}.

Fried described an operation which extends a pseudo-Anosov flow $\Phi_0$ on the complement of a hyperbolic knot $K$ to the closed manifold obtained by $p/q$-surgery along $K$ \cite{Fried83}. The condition for this to apply is given in terms of the degeneracy locus:  

\begin{thm}[Goodman, Fried]
\label{thm: existence of psa flow}
A pseudo-Anosov flow $\Phi_0$ on the complement of $K$ will extend to a pseudo-Anosov flow $\Phi_0(p/q)$ on $K(p/q)$ as long as the absolute value of the algebraic intersection number of $\delta(\Phi_0)$ with the $p/q$ slope of $K$ is at least $2$.
\end{thm}

Gabai and Mosher have independently shown that pseudo-Anosov flows exist on the complement of any hyperbolic link in a closed, orientable 3-manifold. More precisely, they show that given a finite depth taut foliation $\mathcal{F}$ on a compact, connected, orientable, hyperbolic $3$-manifold $M$ with non-empty boundary consisting of tori, there is a pseudo-Anosov flow on the interior of $M$ which is almost transverse to $\mathcal{F}$ (cf. \cite[Theorem C(3)]{Mosher96}).  Unfortunately, no proof has been published, though Landry and Tsang have recently produced the first of several planned articles which will provide a demonstration. See \cite{LT}. 
The following fact about the degeneracy loci of these flows is well known to experts. We include a proof for completeness. 

\begin{thm}[Gabai, Mosher]
\label{thm: degeneracy locus}
Let $K$ be a hyperbolic knot in the $3$-sphere. There is a pseudo-Anosov flow $\Phi_0$ on the complement of $K$ with degeneracy locus $\delta(\Phi_0)$ of the form $b\mu$ or $b \mu + \lambda$ for some integer $b \ne 0$.
\end{thm}

\begin{proof} 
By \cite{Gab83}, there is a finite depth foliation $\mathcal{F}$ on $X(K)$ whose depth zero leaf is a minimal genus Seifert surface for $K$. Then by the results of Gabai and Mosher discussed above, there exists a pseudo-Anosov flow on the interior of $X(K)$, denoted by $\Phi_0$, that is almost transverse to $\mathcal{F}$. Since $K(1/0) \cong S^3$ has trivial fundamental group, it supports no pseudo-Anosov flow. Theorem \ref{thm: existence of psa flow} then implies that $\delta(\Phi_0)$ must be of the form $b\mu$ or $b \mu + \lambda$ for some integer $b$. In the latter case, $b \ne 0$ since the flow is almost transverse to the foliation, in particular the depth zero leaf \footnote{See \cite[Definition 6.45]{Cal07} for the definition of a pseudo-Anosov flow being almost transverse to a foliation. }.
\end{proof}

\begin{thm}
\label{thm: hyp psA 2}
Let $K$ be a hyperbolic knot and $\Phi_0$ a pseudo-Anosov flow on its complement $S^3\setminus K$ with degeneracy locus  $\delta(\Phi_0)$ of the form $b\mu$ or $b \mu + \lambda$ for some integer $b \geq 1$. 
\begin{enumerate}[leftmargin=*]
\setlength\itemsep{0.3em}
\item[{\rm (1)}] If $\delta(\Phi_0) = b \mu$ and $K(1/q)$ or $K(2/q)$ is not $LO$ then $bq = \pm 1$.
\item[{\rm (2)}] If $\delta(\Phi_0) = b\mu + \lambda$ and $K(1/q)$ is not $LO$ then $bq \in \{1, 2\}$.
\item[{\rm (3)}] If $\delta(\Phi_0) = b\mu + \lambda$ and $K(2/q)$ is not $LO$  then $bq \in \{1,2, 3\}$.
\end{enumerate}
\end{thm}

\begin{proof}
For $p = 1$ or $2$, $H_1(K(p/q))$ is a $\mathbb Z/2$ vector space and if $|\delta(\Phi_0) \cdot (p\mu + q \lambda)| \geq 2$ then Theorem \ref{thm: existence of psa flow} implies that $K(p/q)$ admits a pseudo-Anosov flow. The reader will verify that this inequality holds unless $bq$ is as stated in (1), (2) and (3). Corollary \ref{cor: sum of z2s} then shows that $\pi_1(K(p/q))$ is left-orderable, which completes the proof.
\end{proof}

\begin{proof}[Proof of Theorem \ref{thm: $LO$ surgery small p full}]
Suppose that $K$ is a non-trivial knot in $S^3$, $p = 1$ or $2$, and $q \ne 0$, and $K(p/q)$ is not $LO$. As a non-trivial knot, $K$ is either a torus knot, a hyperbolic knot, or a satellite knot. We consider these three cases separately.
    
The known results on left-orderable Dehn surgeries on torus knots stated at the beginning of \S \ref{sec: $LO$ surgery for p small} imply that $K = T(2, 3\varepsilon)$ for some $\varepsilon \in \{\pm 1\}$ and $p/q \in \{\varepsilon, 2\varepsilon\}$. This yields case (1) of the theorem. 

Suppose next that $K$ is hyperbolic. Theorem \ref{thm: existence of psa flow}(2) implies that there is a pseudo-Anosov flow $\Phi_0$ on the complement of $K$ with degeneracy locus $\delta(\Phi_0)$ of the form $b\mu$ or $b \mu + \lambda$ for some integer $b \ne 0$. If $b \geq 1$, Theorem \ref{thm: hyp psA 2} implies that if $p = 1$ then $1/q$ is contained in either $\{-1, 1\}$ or $\{1/2, 1\}$. If $b \leq -1$, replace $K$ by its mirror image $K^*$ and note that $\Phi_0$ is transformed into a pseudo-Anosov flow $\Phi_0^*$ on $S^3 \setminus K^*$ with degeneracy locus $\delta(\Phi_0^*)$ of the form $b^*\mu$ or $b^* \mu + \lambda$ with $b^* = -b \geq 1$. Then as $K^*(-1/q) \cong K(1/q)$, $-1/q$ is contained in either $\{-1, 1\}$ or $\{1/2, 1\}$ and therefore $1/q$ is contained in either $\{-1, 1\}$ or $\{-1, -1/2\}$.

In the case $p = 2$, if $b \ge 1$ then Theorem \ref{thm: hyp psA 2} implies that $2/q$ is contained in either $\{\pm 2\}$ or $\{2/3, 2\}$. Arguing as in the $p = 1$ case above if $b \le -1$ completes the proof in Case (3).

Finally suppose that $K$ is a satellite knot $P(K_0)$ and that $K(p/q)$ is not $LO$. Then $K(p/q) = P(p/q) \cup _{T} X_0$ where $X_0$ is the exterior of $K_0$ and $T = \partial P(p/q) = \partial X_0$. Since $|H_1(K(p/q))| = p$, if $T$ is incompressible in $P(p/q)$ then $K(p/q)$ would be $LO$ by Theorem \ref{thm: $LO$ small HRRW15}, contrary to our hypotheses.

Thus $T$ compresses in $P(p/q)$. Then by Theorem \ref{thm: sch red} either (i) $P = C_{m,n}$ and $p/q$ is the cabling slope $mn$, or (ii) $P(p/q)$ is a solid torus.

If possibility (i) arises, then $q = \varepsilon \in \{\pm 1\}$ and therefore $mn = \varepsilon p$. Since $m \geq 2$, we must have $p = m = 2$ and $n = \varepsilon$. Conversely, when $p, q, m, n$ take on these values, $K(p/q)$ has an $\mathbb RP^3$ summand, so is not $LO$. To complete the proof of Case (2) we must rule out possibility (ii).

Assume that (ii) holds. Then Theorems \ref{thm: sch red} and \ref{thm: C(T)} imply that either (a) $P = C_{m,n}$ and $\Delta(p/q, mn/1) = 1$, or (b) $P$ is a 1-bridge braid and $|q| = 1$. It is easy to verify that the only solutions to (a) are $p = 1, q = \varepsilon, m = 2, n = \varepsilon$, and $p = 2, q = \varepsilon, m = 3, n = \varepsilon$, for some $\varepsilon = \pm 1$. Thus in both cases (a) and (b), $|q| = 1$. Moreover, in case (a) $w = w(P) = m \geq 2$, and in case (b) $w \ge 4$ by Lemma \ref{lemma: 1 bridge winding number}.

We have $K(p/q) = K(p) = K_0(p/w^2)$ by Theorem \ref{thm: C(T)}, where $w^2 \geq 4$. If $K_0$ is a torus knot then $K_0(p/w^2)$ is $LO$ since $p/w^2 < 1$. The same conclusion holds if $K_0$ is hyperbolic by Case (3). Finally, the conclusion holds if $K_0$ is a satellite knot by the argument above that any non-integral surgery $p/q$ on a satellite knot with $p = 1$ or $2$ is $LO$. 
\end{proof}

\begin{cor}
    \label{cor: LO fibered p=1 2}
If $K$ is a non-trivial fibred knot, then $K(1/q)$ and $K(2/q)$ are $LO$ for $|q| \geq 2$. 
\end{cor}
\begin{proof}
If $K$ is hyperbolic, then the monodromy of $K$ is pseudo-Anosov (\cite{Thurston98}) and hence its suspension flow gives a pseudo-Anosov flow on the complement of $K$. Up to replacing $K$ by its mirror image the degeneracy locus of the suspension flow  is either $b \mu$ for $b \geq 1$ or $b \mu + \lambda$ for some integer $|b| \geq 2$ by \cite[Theorem 8.8]{Gab97}. By Theorem \ref{thm: hyp psA 2}, $K(1/q)$ and $K(2/q)$ are $LO$ for $|q| \geq 2$. 

The same conclusion holds when $K$ is either a torus knot or satellite knot by Theorem \ref{thm: $LO$ surgery small p full}(1) and (2) respectively.
\end{proof}

By \cite[Corollary 6.17]{BGH2}, if $K$ is a knot whose complement admits a pseudo-Anosov flow with non-meridional degeneracy locus then all the $n$-fold cyclic branched covers of $K$, $n \ge 2$, are $LO$. Theorem \ref{thm: hyp psA 2} therefore gives the following.

\begin{cor}
    \label{cor: cyclic cover}
Let $K$ be a hyperbolic knot such that some $n$-fold cyclic branched cover of $K$, $n \ge 2$, is not $LO$. Then $K(1/q)$ and $K(2/q)$ are $LO$ for $|q| \geq 2$. 
\end{cor} 

Since the double branched cover of an alternating knot is not $LO$ \cite{BGW13} we obtain
 
\begin{cor} 
    \label{cor: LO alternating p=1 2}
If $K$ is a non-trivial alternating knot, then $K(1/q)$ and $K(2/q)$ are $LO$ for $|q| \geq 2$. 
\end{cor} 
\begin{proof}
By \cite{Menasco84}, an alternating knot cannot be a satellite knot. 
The conclusion follows immediately from Corollary \ref{cor: cyclic cover} if $K$ is hyperbolic. If $K$ is a torus knot, the conclusion follows from Theorem \ref{thm: $LO$ surgery small p full}(1).
\end{proof}

\begin{remark}
Note that Corollaries \ref{cor: LO fibered p=1 2} and \ref{cor: LO alternating p=1 2} can also be proven using the universal circle actions (see \cite{CD03, BH19})  associated to the taut foliations on the surgered manifolds constructed by Roberts in \cite{RobertsSurfacebundle2} and \cite{Roberts1995} respectively.  The Euler class of these universal circle actions vanishes by \cite[Proposition 2.1]{Hu2} and therefore the fundamental groups of the surgered manifolds are LO (cf. \cite[Theorem 8.3]{BH19}).
\end{remark}

For general $p/q \in \mathbb Q$, we can find a pseudo-Anosov flow $\Phi(p/q)$ on $K(p/q)$ as long as $|\delta(\Phi) \cdot (p\mu + q \lambda)| \geq 2$ and consequently faithful representations $\rho_{p/q}: \pi_1(K(p/q)) \to \mbox{Homeo}_+(S^1)$ (see the Appendix). As noted in the Appendix, $\pi_1(K(p/q))$ is left-orderable if the Euler class $e(\rho_{p/q}) \in H^2(K(p/q))$ vanishes. Equivalently, the Euler class of the normal bundle to the flow $e(\nu_{\Phi(p/q)})$ vanishes (Proposition \ref{prop: equal euler classes}).  But this does not appear to happen very frequently. For instance, though the methods of \cite{Hu2} can be often be used to show that there are infinitely many values of $p/q$ for which $e(\rho_{p/q}) = 0$, they also show that the set of such $p/q$ is typically nowhere dense in the reals.

\appendix
\section{The Euler class of Fenley's asymptotic circle}
\label{app: fenley}
The goal of this appendix is to calculate the Euler class of Fenley's asymptotic circle representation in terms of the topology of the associated pseudo-Anosov flow (Proposition \ref{prop: equal euler classes}). The result is presumably known to experts, though we do not know a reference.

Given a pseudo-Anosov flow $\Phi$ on a closed, connected, orientable $3$-manifold $W$, let $\widetilde \Phi$ be the pull-back of $\Phi$ to the universal cover $\widetilde W$ of $W$.

\begin{thm}
{\rm (\cite[Proposition 4.2]{FM01})}
\label{thm: orbit space}
The orbit space $\mathcal{O}$ of $\widetilde \Phi$ is homeomorphic to $\mathbb R^2$. Moreover, the projection $\pi: \widetilde W \to \mathcal{O}$ is a locally-trivial fibre bundle whose flow line fibres are homeomorphic to $\mathbb R$.
\end{thm}
The action of $\pi_1(W)$ on $\widetilde W$ descends to one on $\mathcal{O}$ by homeomorphisms. Since the flow lines in $\widetilde W$ inherit a coherent $\pi_1(W)$-invariant orientation, the action of $\pi_1(W)$ on $\mathcal{O}$ is by orientation-preserving homeomorphisms, so we obtain a homomorphism
$$\psi: \pi_1(W) \to \mbox{Homeo}_+(\mathcal{O})$$
Fenley has constructed an ideal boundary for $\mathcal{O}$ over which this action extends.

\begin{thm}
{\rm (\cite[Theorem A]{Fen12})}
\label{thm: fenley's univ circle}
There is a natural compactification $\mathcal{D} = \mathcal{O} \cup \partial \mathcal{O}$ of $\mathcal{O}$ where $\mathcal{D}$ is homeomorphic to a disk with boundary circle $\partial \mathcal{O}$. The action of $\pi_1(W)$ on $\mathcal{O}$ extends to one on $\mathcal{D}$ by homeomorphisms.
\end{thm}
It follows from Fenley's construction that the action of $\pi_1(W)$ on the ideal boundary $\partial \mathcal{O}$ of $\mathcal{O}$ is faithful. That is, the associated homomorphism
$$\rho_\Phi: \pi_1(W) \to \mbox{Homeo}_+(\partial \mathcal{O})$$
is injective. We think of $\rho_\Phi$ as taking values in $\mbox{Homeo}_+(S^1)$.

The action of $\pi_1(W)$ on $\mathcal{O}$ gives rise to a diagonal action on $\widetilde W \times \mathcal{O}$ via
$$\gamma \cdot (x, y) = (\gamma \cdot x, \psi(\gamma)(y)),$$
which is equivariant with respect to the projection $\widetilde W \times \mathcal{O} \to \widetilde W$. Taking quotients determines a
locally-trivial $\mathcal{O}$-bundle
$$W \times_{\psi} \mathcal{O} = \big(\widetilde W \times \mathcal{O} \big)/ \pi_1(W) \to W$$

\begin{lemma}
\label{lemma: equivalence}
The $\mathbb R^2$-bundle $W \times_{\psi} \mathcal{O}\to W$ is topologically equivalent to the normal bundle $\nu(\Phi)$ of $\Phi$.
\end{lemma}

\begin{proof}
We use $E(\nu(\Phi))$ and $E(\nu(\widetilde \Phi))$ to denote the total spaces of $\nu(\Phi)$ and $\nu(\widetilde \Phi)$.

Fix a Riemannian metric on $W$ and pull it back to $\widetilde W$. Since the group of deck transformations of the cover $\widetilde W \to W$ acts as isometries of this metric, the associated exponential map $\exp: T \widetilde W \to \widetilde W$ is $\pi_1(W)$-equivariant, as is the composition
$$\omega: E(\nu(\widetilde \Phi)) \xrightarrow{\exp} \widetilde W \xrightarrow{\pi} \mathcal{O}$$
Similarly if $E_\epsilon(\nu(\widetilde \Phi))$ denotes the set of vectors of length less than $\epsilon$ in $E(\nu(\widetilde \Phi))$, then $E_\epsilon(\nu(\widetilde \Phi))$ is $\pi_1(W)$-invariant and the restriction
$$\omega_\epsilon : E_\epsilon(\nu(\widetilde \Phi)) \xrightarrow{\omega|_{E_\epsilon(\nu(\widetilde \Phi))}} \mathcal{O}$$
is $\pi_1(W)$-equivariant. Since $W$ is compact and $\widetilde \Phi$ is invariant under the action of $\pi_1(W)$, there is an $\epsilon > 0$ such that for each $x \in \widetilde W$,
$\omega_\epsilon$ determines a homeomorphism between the open disk fibre of $E_\epsilon(\nu(\widetilde \Phi))$ over $x$ and its image in $\mathcal{O}$.

Fix a diffeomorphism $\varphi: ([0, \infty), 0) \to ([0, \epsilon), 0)$ and consider the equivariant embedding
$$E(\nu(\widetilde \Phi)) \xrightarrow{(x, v) \mapsto (x, \varphi(\|v\|)v)} E_\epsilon(\nu(\widetilde \Phi)) \xrightarrow{(x, v) \mapsto (x, \omega_\epsilon(v))} \widetilde W \times \mathcal{O}$$
which sends fibres into open subsets of fibres. Quotienting out by $\pi_1(W)$ yields an embedding
$$E(\nu(\Phi)) \to  W \times_{\psi} \mathcal{O}$$
which also sends fibres into open subsets of fibres. Let $E_\nu$ be the image of this embedding and $W_0 \subset  E_\nu \subset W \times_{\psi} \mathcal{O}$ the image of the zero section of $E(\nu(\Phi))$. Since $W_0 \to E_\nu \to W_0$ is a sub-$\mathbb R^2$-microbundle of the $\mathbb R^2$-microbundle $W_0 \to W \times_{\psi} \mathcal{O} \to W_0$, the main result of \cite{Kister64} implies that $E_\nu \to W_0$ and $W \times_{\psi} \mathcal{O} \to W_0$ are isomorphic
$\mathbb R^2$-bundles, which completes the proof.
\end{proof}

\begin{prop}
\label{prop: equal euler classes}
$e(\rho_\Phi) =  e(\nu(\Phi))$
\end{prop}

\begin{proof}
The Euler class of $\rho_\Phi$ coincides with that of the oriented circle bundle
$$W \times_{\rho_\Phi} \partial \mathcal{O} = \big(\widetilde W \times \partial \mathcal{O} \big)/ \pi_1(W) \to W$$
by \cite[Lemma 2]{Mil58}, while  the Euler class of this circle bundle coincides with that of the associated $2$-disk bundle $\big(\widetilde W \times \mathcal{D} \big)/ \pi_1(W) \to W$ by \cite[\S 5.7]{Spa66} and therefore to that of $W \times_{\psi} \mathcal{O} \to W$. Lemma \ref{lemma: equivalence} then
completes the proof.
\end{proof}

\begin{cor}
\label{cor: pseudo-Anosov and $LO$ 2}
If $W$ is a closed, connected, orientable $3$-manifold which admits a pseudo-Anosov flow $\Phi$ whose normal bundle has zero Euler class, then
$\pi_1(W)$ is left-orderable.
\end{cor}

\begin{proof}
Let $\rho_\Phi: \pi_1(W) \to \mbox{Homeo}_+(S^1)$ be Fenley's universal circle representation and note that as $e(\rho_\Phi) = e(\nu(\Phi)) = 0$, $\rho_\Phi$ lifts to a representation $\widetilde \rho: \pi_1(W) \to \mbox{Homeo}_{\mathbb Z}(\mathbb R) \leq \hbox{Homeo}_+(\mathbb R)$. Since $\rho_\Phi$ is faithful, so is $\widetilde \rho$. Finally, as a closed manifold admitting a pseudo-Anosov flow, $W$ is irreducible and therefore $\pi_1(W)$ is left-orderable by \cite[Theorem 1.1]{BRW05}.
\end{proof}

\begin{cor}
\label{cor: sum of z2s}
Let $W$ be a rational homology $3$-sphere for which $H_1(W)$ is a $\mathbb Z/2$ vector space. If $W$ admits a pseudo-Anosov flow $\Phi$, then $\pi_1(W)$ is left-orderable.
\end{cor}

\begin{proof}
By Corollary \ref{cor: pseudo-Anosov and $LO$ 2}, we know that $\pi_1(W)$ will be left-orderable as long as $e(\nu(\Phi)) = 0$. This is obvious if $H^2(W) \cong H_1(W) \cong \{0\}$ and holds if $H_1(W)$ is a $\mathbb Z/2$ vector space by \cite[Proposition 2.1]{Hu2}.
\end{proof}

{
\footnotesize
\bibliographystyle{abbrvnat}
\bibliography{bgh_lo}
}

\end{document}